\documentclass[a4paper,fleqn]{cas-dc}

\usepackage[authoryear]{natbib}

\usepackage{lipsum}
\usepackage{tikz}

\usepackage{graphicx}
\usepackage{caption}
\usepackage{subcaption}

\usepackage{float}
\usepackage[siunitx]{circuitikz}
\usepackage{placeins}
\usepackage{needspace}

\def\tsc#1{\csdef{#1}{\textsc{\lowercase{#1}}\xspace}}
\tsc{WGM}
\tsc{QE}

\begin{document}
\let\WriteBookmarks\relax
\def\floatpagepagefraction{1}
\def\textpagefraction{.001}

\shorttitle{3D PINN for Pulmonary Hemodynamics using RCR}
\shortauthors{C. Taboada et~al.}

\title [mode = title]{A 3D PINN-Based Framework for Patient-Specific Pulmonary Hemodynamics Using Windkessel Models}                      



\author[1]{Celia Taboada}[orcid=0009-0004-8170-7944]
\cormark[1]
\ead{celia.taboada@upm.es}
\credit{Conceptualization, Methodology, Software, Formal analysis, Investigation, Writing - original draft, Visualization}

\affiliation[1]{organization={Department of Continuum Mechanics and Theory of Structures},
                addressline={Prof. Aranguren 3}, 
                city={Madrid},
                postcode={28040}, 
                country={Spain}}

\author[1]{Luis Nuñez}
\credit{Investigation, Software, Validation}

\author[2]{Angélica Casero}[orcid=0009-0005-4912-7568]
\ead{angelica.caserobaeza@unibe.ch}
\credit{Conceptualization, Supervision, Data curation, Writing - review \& editing. }

\affiliation[2]{organization={ARTORG Center for Biomedical Engineering Research, University of Bern},
                addressline={Freiburgstrasse 3}, 
                city={Bern},
                postcode={3010}, 
                country={Switzerland}}

\author[1]{Pedro Navas}[orcid=0000-0002-0482-378X]
\ead{pedro.navas@upm.es}
\ead[URL]{https://www1.caminos.upm.es/estructuras2/en/}
\credit{Conceptualization, Supervision, Funding acquisition, Resources, Project administration, Writing - review \& editing}

\cortext[cor1]{Corresponding author}

\begin{abstract}
Cardiovascular fluid dynamics modeling relies heavily on mesh-based numerical methods such as the Computational Fluid Dynamics (CFD). While accurate, these methods entail significant computational costs and complex meshing workflows that limit their translation to clinical practice. Physics-Informed Neural Networks (PINNs) are a promising meshless alternative that integrates the governing equations directly into the training of a neural network. In this work, we propose a PINN framework for pulsatile blood flow in patient-specific geometries with dynamic lumped-parameter Windkessel (RCR) outlet conditions, which represent the distal vascular impedance.
 
Coupling the Windkessel ordinary differential equation to a 3D network raises two difficulties: the ODE alone does not fix the pressure level, and the local pressure variations that drive the flow are orders of magnitude smaller than the basal pressure imposed by the outlet. We address them with a decomposition of the pressure into a time-dependent basal part and a spatially varying part, an initial condition and gauge for the Windkessel equation, a global mass-balance constraint, and a periodic Fourier embedding of time. The framework is validated on two-dimensional canonical channels (steady and pulsatile, with and without RCR outlet) and then applied to a patient-specific right pulmonary artery, where it reproduces the outlet flow and mean pressure, and the velocity and dynamic pressure fields, with errors of the order of the data resolution.
 
Finally, we solve the inverse problem of identifying the three Windkessel parameters of the outlet from synthetic data, recovering them with errors of about $0.3\%$ from deliberately wrong initial values. Residual biases in the local pressure field and the need to test robustness to noisy clinical data remain open, but the results support PINNs as complementary tools for meshless hemodynamics and for patient-specific digital twins.

\end{abstract}

\begin{graphicalabstract}
\includegraphics[width=\linewidth]{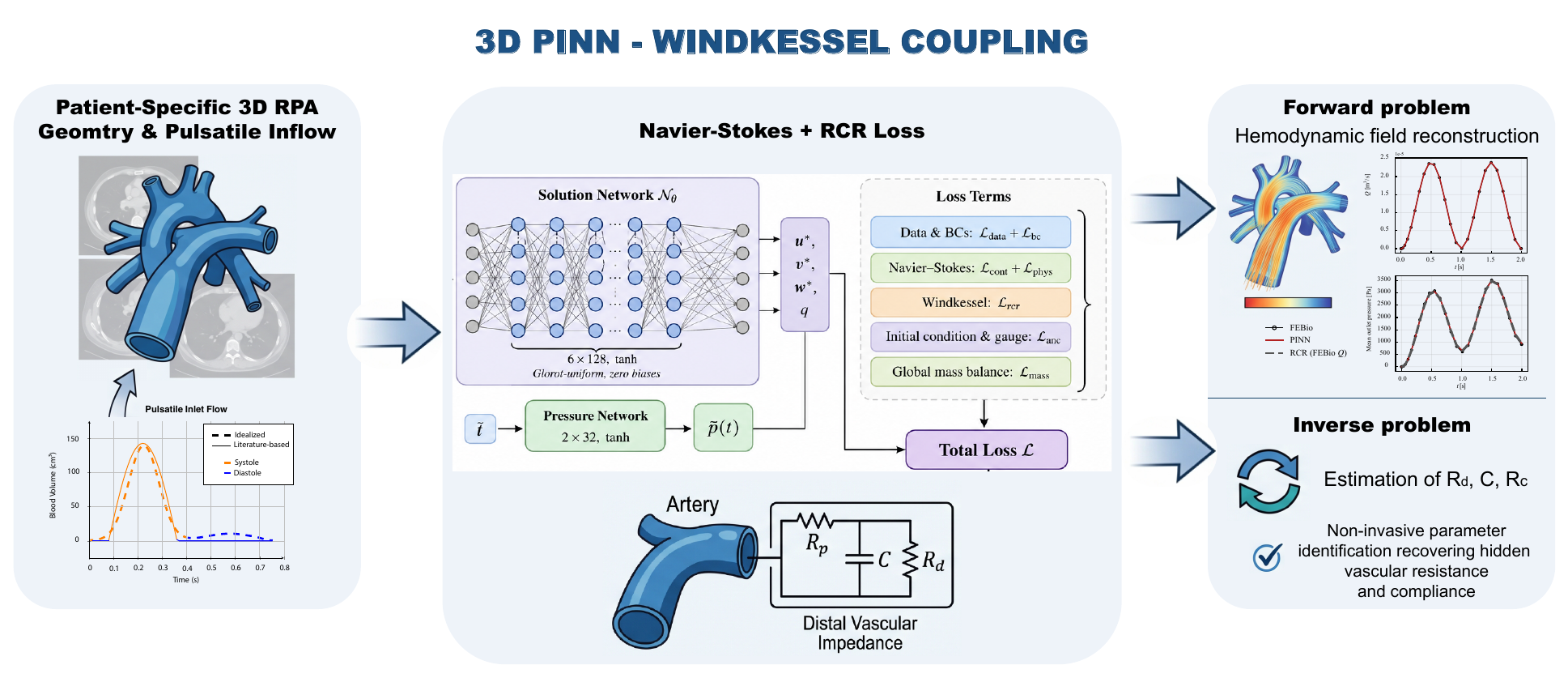}
\end{graphicalabstract}

\begin{highlights}
\item PINN framework developed for 3D pulsatile hemodynamics in pulmonary arteries.
\item RCR Windkessel boundary conditions integrated into neural network training.
\item Addresses pressure-scale challenges in 3D lumped-parameter PINN coupling.
\end{highlights}


\maketitle











\section{Introduction}\label{sec:intro}
 
Cardiovascular diseases remain a leading cause of mortality worldwide, and pulmonary vascular disease is a particularly demanding case. Pulmonary hypertension, for instance, is defined by an elevated mean pulmonary arterial pressure, and its definitive diagnosis and follow-up still rely on invasive right-heart catheterization \cite{humbert2022guidelines}. Implantable pressure sensors such as CardioMEMS\texttrademark\ have opened the door to monitoring \cite{abraham2017implantable}, but they provide pressure at a single location and say nothing about the underlying flow field. Patient-specific hemodynamic modeling offers a complementary, non-invasive route: starting from medical imaging and a few clinical measurements, it can estimate the full pressure and velocity fields, as well as derived indices such as wall shear stress or pressure drops, that cannot be measured directly \cite{casero2026patient}.
 
The reference tool for this task is image-based computational fluid dynamics (CFD) or, equivalently, finite element analysis of the Navier--Stokes equations \cite{taylor2009patient}. These mesh-based solvers are accurate and well validated, but their clinical translation is hindered by several factors. Each new patient requires segmentation, surface repair and volumetric meshing, which are labor-intensive and difficult to automate. Each simulation must be repeated from scratch for every change in geometry, inflow or boundary condition. Moreover, the outcome is strongly conditioned by the outflow boundary conditions, which must be calibrated for every case through lumped-parameter (Windkessel) models \cite{vignonclementel2006outflow, westerhof2009arterial}. Finally, classical solvers are ill suited to assimilate noisy or sparse measurements, such as those provided by 4D flow MRI or Doppler ultrasound, and they cannot solve inverse problems without wrapping them in costly optimization loops.
 
Physics-Informed Neural Networks (PINNs) have emerged as an alternative that addresses several of these limitations \cite{raissi2019physics, karniadakis2021physics}. A PINN represents the solution fields with a neural network whose inputs are the space-time coordinates, and it is trained by minimizing a loss that contains the residual of the governing equations, evaluated through automatic differentiation at scattered collocation points, together with boundary conditions and, optionally, measured data. The method is meshless, treats forward, inverse and data-assimilation problems within the same formulation, and yields a continuous, differentiable representation of the solution. These advantages come at a price: PINNs are notoriously difficult to train, mainly because of the spectral bias of fully connected networks \cite{rahaman2019spectral, tancik2020fourier}, the competition between loss terms of very different magnitude \cite{wang2021understanding}, and the stiffness of the resulting optimization problem \cite{krishnapriyan2021characterizing}. All of them are exacerbated in pulsatile cardiovascular flows, where time-periodic dynamics, convective and viscous terms of different scales, and a pressure that is defined only up to a level set by the boundary conditions coexist.
 

Focusing specifically on cardiovascular fluid mechanics, PINNs have been progressively tailored to reconstruct complex arterial flow patterns, estimate near-wall hemodynamic metrics, and perform super-resolution or denoising on clinical imaging data \cite{arzani2021physics, kissas2020machine, fathi2022super, arzani2023uncovering}. Recent studies have demonstrated their efficacy in capturing 3D spatio-temporal velocity and pressure fields across complex arterial bifurcations and patient-specific topologies \cite{sarma2024physics, gonzalez2025reconstruction, aghaee2025pinning}. Furthermore, growing attention has been directed toward integrating reduced-order or lumped-parameter models into PINN architectures to capture downstream vascular resistance and compliance \cite{ames2024coupling, muller2025uncertainty, vergara2024coupling, park2024estimation}. However, enforcing dynamic differential Windkessel boundary conditions within 3D patient-specific domains introduces severe numerical stiffness and pressure-velocity scale imbalances that remain insufficiently addressed in the current literature \cite{albocher2024physics, santos2024evaluation}.
 
The difficulty is structural. The Windkessel relation is an ordinary differential equation in time for the outlet pressure, so it admits a whole family of solutions that differ in the initial condition and in the pressure level; a PINN that only penalizes its residual may satisfy it with a physically wrong pressure and a spuriously small outflow. In addition, the local pressure variations that drive the flow (from one to a few tens of pascals in the problems considered here) are small compared with both the dynamic pressure $\rho U^2$ and the absolute basal pressure imposed by the Windkessel, so a naive non-dimensionalization leaves the pressure field poorly conditioned.


To address these specific optimization hurdles and bridge the gap between physics-informed learning and dynamic boundary condition coupling, this work presents a comprehensive 3D PINN framework for patient-specific pulmonary hemodynamics. Building upon our recent findings on right pulmonary artery modeling \cite{casero2026patient}---which established that rigid vessel walls combined with compliant 3-element Windkessel (3WK) boundary conditions yield high physiological fidelity at minimal computational cost---we embed the dynamic RCR differential equations directly into the network optimization graph. The primary contributions of this paper are threefold: (i) establishing a progressive complexity workflow transitioning from canonical 2D flows to a 3D patient-specific right pulmonary artery domain; (ii) providing a critical quantitative evaluation of the pressure-scale disparity and optimization biases induced by dynamic RCR boundary conditions in 3D PINNs; and (iii) demonstrating the framework's capability to solve inverse data assimilation problems by non-invasively identifying hidden distal vascular resistance and compliance parameters from sparse flow measurements. The overall workflow of the proposed framework is summarized in Fig.~\ref{fig:intro}.

\begin{center}
    \includegraphics[width=\linewidth]{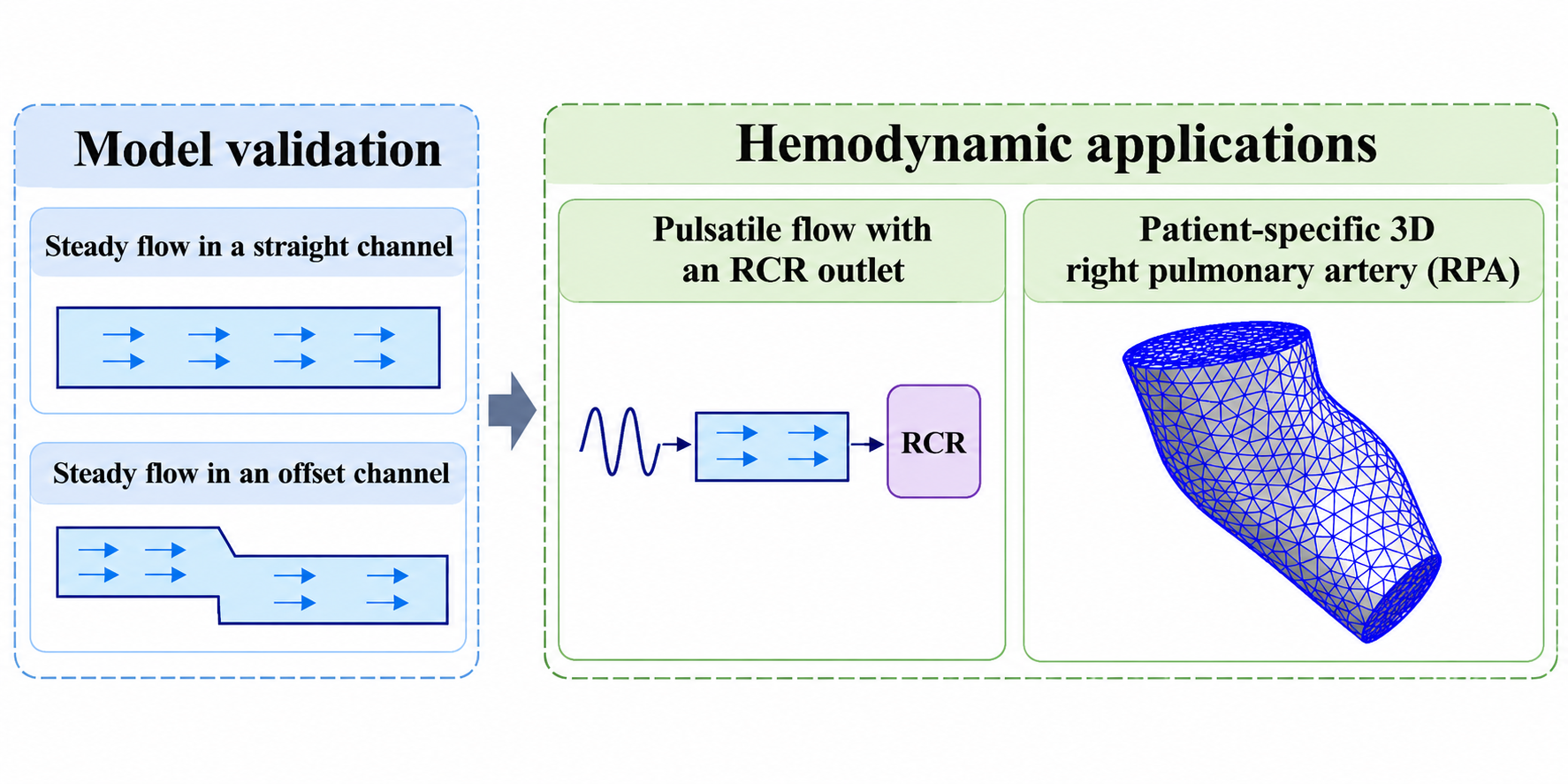}
    \captionof{figure}{Schematic overview of the framework, comprising two stages: (1) model validation under steady-flow conditions in straight and offset channels, and (2) hemodynamic applications involving pulsatile flow with a three-element Windkessel (RCR) outlet boundary condition and a patient-specific three-dimensional model of the right pulmonary artery (RPA).}
    \label{fig:intro}
\end{center}


The remainder of this paper is organized as follows. Section \ref{sec:methods} details the mathematical formulation of blood flow and the RCR Windkessel boundary conditions. Section \ref{subsec:pinns} presents the proposed PINN framework, including the network architecture, the pressure-decomposition non-dimensionalization strategy, and the composite loss function. Section \ref{sec:validation} provides a systematic validation on canonical 2D benchmark cases of increasing complexity: a steady straight channel, a steady offset channel. Section \ref{sec:RPA} then applies the framework to a pulsatile channel coupled with an RCR and an outlet patient-specific 3D right pulmonary artery, addressing first the forward hemodynamic field reconstruction under dynamic RCR coupling and subsequently for the later the inverse identification of hidden vascular parameters. Finally, Section \ref{sec:conclusions} summarizes the main conclusions and outlines future clinical perspectives.

\clearpage
\section{Mathematical and Computational Framework}\label{sec:methods}

\subsection{Governing Equations}\label{subsec:navier_stokes}

Blood flow through large-to-medium vessels, such as the main and right pulmonary arteries (RPA), is characterized as an incompressible, viscous Newtonian fluid governed by the three-dimensional, unsteady Navier-Stokes equations \cite{milnor1989hemodynamics}. Over a spatio-temporal domain $\Omega \subset \mathbb{R}^3 \times [0, T]$ bounded by $\partial\Omega = \Gamma_{in} \cup \Gamma_{wall} \cup \Gamma_{out}$, conservation of mass and linear momentum are expressed as:
\begin{equation}
    \nabla \cdot \mathbf{u} = 0 \quad \text{in } \Omega \times [0, T]
\end{equation}
\begin{equation}
    \rho \left( \frac{\partial \mathbf{u}}{\partial t} + (\mathbf{u} \cdot \nabla)\mathbf{u} \right) = -\nabla p + \mu \nabla^2 \mathbf{u} + \rho \mathbf{b} \quad \text{in } \Omega \times [0, T]
\end{equation}
where $\mathbf{u}(\mathbf{x}, t) = (u, v, w)^T$ is the fluid velocity vector field, $p(\mathbf{x}, t)$ is the static pressure, $\rho = 1060 \text{ kg/m}^3$ denotes the blood mass density, $\mu = 0.004 \text{ Pa}\cdot\text{s}$ represents the dynamic viscosity, and $\mathbf{b}$ corresponds to body forces, which are neglected in this study \cite{milnor1989hemodynamics}.

The physical regime of the blood flow within the pulmonary artery is dictated by two key non-dimensional parameters:
\begin{itemize}
    \item \textbf{Reynolds Number ($Re$):} Evaluated as $Re = \frac{\rho v D}{\mu}$, where $D$ is the hydraulic vessel diameter. Based on patient-specific cardiac output ($CO = 6.183 \times 10^{-5} \text{ m}^3/\text{s}$), the mean flow velocity $v \approx 0.061 \text{ m/s}$ yields a mean Reynolds number of $Re \approx 429$ \cite{casero2026patient}. Even under conservative peak systolic flow conditions ($140 \text{ cm}^3/\text{s}$), $Re$ remains below $1770$ \cite{casero2026patient}. This confirms that, despite potential localized secondary flows induced by arterial curvature, the flow regime remains predominantly laminar \cite{casero2026patient}.
    \item \textbf{Womersley Number ($\alpha$):} Calculated as $\alpha = r \sqrt{\frac{\omega \rho}{\mu}}$, where $r$ is the vessel radius and $\omega = 2\pi \frac{HR}{60}$ is the angular frequency associated with the heart rate ($HR = 80 \text{ bpm}$) \cite{casero2026patient}. For the analyzed right pulmonary artery ($r = 1.33 \text{ cm}$), the resulting Womersley number is $\alpha \approx 8.33$ \cite{casero2026patient}. A value of $\alpha > 8$ confirms that transient, pulsatile inertia strongly dominates over viscous diffusion, ruling out quasi-steady flow assumptions and requiring an explicitly unsteady temporal treatment \cite{casero2026patient}.
\end{itemize}

\begin{figure}
  \centering
    \includegraphics[width=\linewidth]{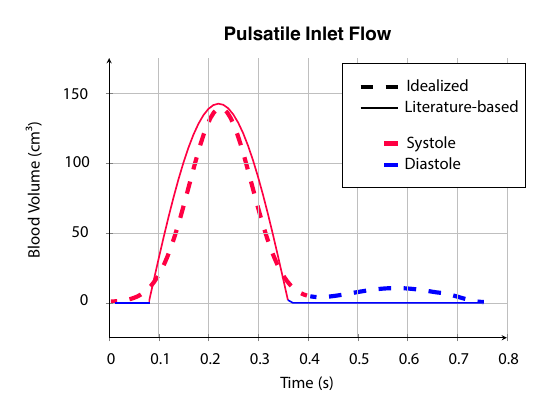}
    \caption{Comparison of inflow boundary condition waveforms prescribed at $\Gamma_{in}$: patient-adapted physiological profile derived from clinical literature versus the simplified half-sine pulsatile curve. Adapted from \cite{casero2026patient}.}\label{fig:inflow_waveform}
\end{figure}

To complete the mathematical boundary value problem, specific constraints are assigned across the domain boundaries:
\begin{enumerate}
    \item \textbf{Inlet ($\Gamma_{in}$):} A time-dependent, pulsatile volumetric inflow profile $Q_{in}(t)$ is enforced. As illustrated in Fig. \ref{fig:inflow_waveform}, the inlet flow can be modeled using a literature-based physiological profile or an idealized half-sine function scaled to match the total stroke volume ($SV = CO / HR$) during the systolic period ($T_{sys} \approx 0.375 T$) \cite{casero2026patient}.
    \item \textbf{Vessel Wall ($\Gamma_{wall}$):} A non-slip Dirichlet boundary condition is prescribed:
    \begin{equation}
        \mathbf{u}(\mathbf{x}, t) = \mathbf{0} \quad \text{on } \Gamma_{wall} \times [0, T]
    \end{equation}
    While Fluid-Structure Interaction (FSI) captures elastic vessel wall mechanics, assuming rigid walls significantly mitigates computational complexity. Recent comparative validations against clinical sensor data demonstrate that assuming a rigid wall introduces a negligible discrepancy in peak pressure ($\approx 1.1\%$, corresponding to a $0.4 \text{ mmHg}$ offset) relative to full FSI solutions, validating its use in computationally efficient frameworks \cite{casero2026patient}.
\end{enumerate}

\subsection{Windkessel (RCR) Boundary Conditions}\label{subsec:rcr}

Truncating the 3D computational domain at distal boundaries without imposing artificial pressure reflections requires phantom boundary conditions that capture downstream vascular impedance \cite{westerhof2009arterial}. To achieve physiological realistic pressure waveforms at the outlets $\Gamma_{out}$, a three-element Windkessel (3WK / RCR) lumped-parameter model is coupled to the fluid equations \cite{frank1899grundform, westerhof2009arterial}.

\begin{figure}
    \centering
    \subfloat[Cardiovascular system as a closed hydraulic circuit. Adapted from \cite{casero2026patient}.]{%
      \includegraphics[width=0.95\linewidth]{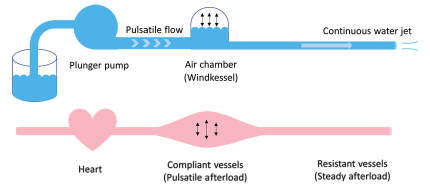}%
        \label{fig:dibujo2}%
   }
    
    \vspace{0.2cm}
    
    \subfloat[Circuit model.]{%
        \begin{circuitikz}[scale=0.68, transform shape]
            \draw (0,0) to[open] (0,2) 
           to[R, l=$R_p$] (2,2) 
           to[short] (2.5,2)
           to[short] (3,2)
            to[short] (3,1.5)
            to[C, l=$C$] (3,0) 
            (3,2) to[short] (5,2)
            to[R, l=$R_d$] (5,0) 
           to[short] (0,0); 
        \end{circuitikz}%
        \label{fig:circuito}%
    }
    \caption{Schematic representation of the 3-element Windkessel (RCR) lumped-parameter circuit coupled at the outlet boundary $\Gamma_{out}$ to model proximal resistance ($R_p$), distal compliance ($C$), and distal vascular bed resistance ($R_d$). Adapted from \cite{casero2026patient, westerhof2009arterial}.}
    \label{fig:rcr_circuit}
\end{figure}

As depicted in Fig. \ref{fig:rcr_circuit}, the RCR circuit models the downstream systemic or pulmonary bed using three hydraulic components: a proximal characteristic resistance $R_p$, an arterial compliance $C$, and a distal peripheral resistance $R_d$ \cite{westerhof2009arterial}.

The differential relationship governing the instantaneous outlet pressure $P_{out}(t)$ and the spatial outlet flow rate $Q_{out}(t) = \int_{\Gamma_{out}} \mathbf{u} \cdot \mathbf{n} \, \mathrm{d}\Gamma$ is derived from circuit conservation laws as:
\begin{equation}\label{eq:rcr_ode}
    C \frac{\mathrm{d}P_{out}(t)}{\mathrm{d}t} + \frac{P_{out}(t)}{R_d} = Q_{out}(t) \left( 1 + \frac{R_p}{R_d} \right) + R_p C \frac{\mathrm{d}Q_{out}(t)}{\mathrm{d}t}
\end{equation}

The total vascular resistance $R_T = R_p + R_d$ is calibrated using the patient's mean pulmonary pressure $P_{mean}$ and mean outlet discharge $Q_{RPA}$ via $R_T = P_{mean} / Q_{RPA}$ \cite{casero2026patient}. In the patient-specific calibration of \cite{casero2026patient}, $R_T$ is partitioned into a proximal resistance $R_p$ ($5\%$ of $R_T$), representing local vessel impedance, and a distal microvascular resistance $R_d$ ($95\%$ of $R_T$), in order to minimize high-frequency pressure-wave reflections at the computational boundary. In the benchmarks of this paper the parameters are instead fixed to the values of the corresponding FE models (Table~\ref{tbl:rcr_2d} and Section~\ref{subsec:3d_inverse}), so that the PINN and the reference solution share exactly the same outlet condition.

The presence of the arterial compliance term $C$ in Eq. \eqref{eq:rcr_ode} is non-negotiable for physiological realism \cite{casero2026patient, ames2024coupling}. When compliance is neglected ($C = 0$), reducing the model to a single-element pure resistance outlet (1WK), the lack of distal volume storage causes extreme, non-physiological pressure spikes (exceeding $130 \text{ mmHg}$ during systole) and unrealistic diastolic drops to $0 \text{ mmHg}$ \cite{casero2026patient}. Conversely, incorporating the full dynamic RCR system captures both the systolic pressure accumulation and the smooth exponential diastolic decay recorded by remote implantable pressure sensors such as CardioMEMS\texttrademark\ \cite{abraham2017implantable, casero2026patient}.

\FloatBarrier
\section{Physics-Informed Neural Networks (PINNs)}\label{subsec:pinns}

\begin{figure*}
    \centering
    \includegraphics[width=\linewidth]{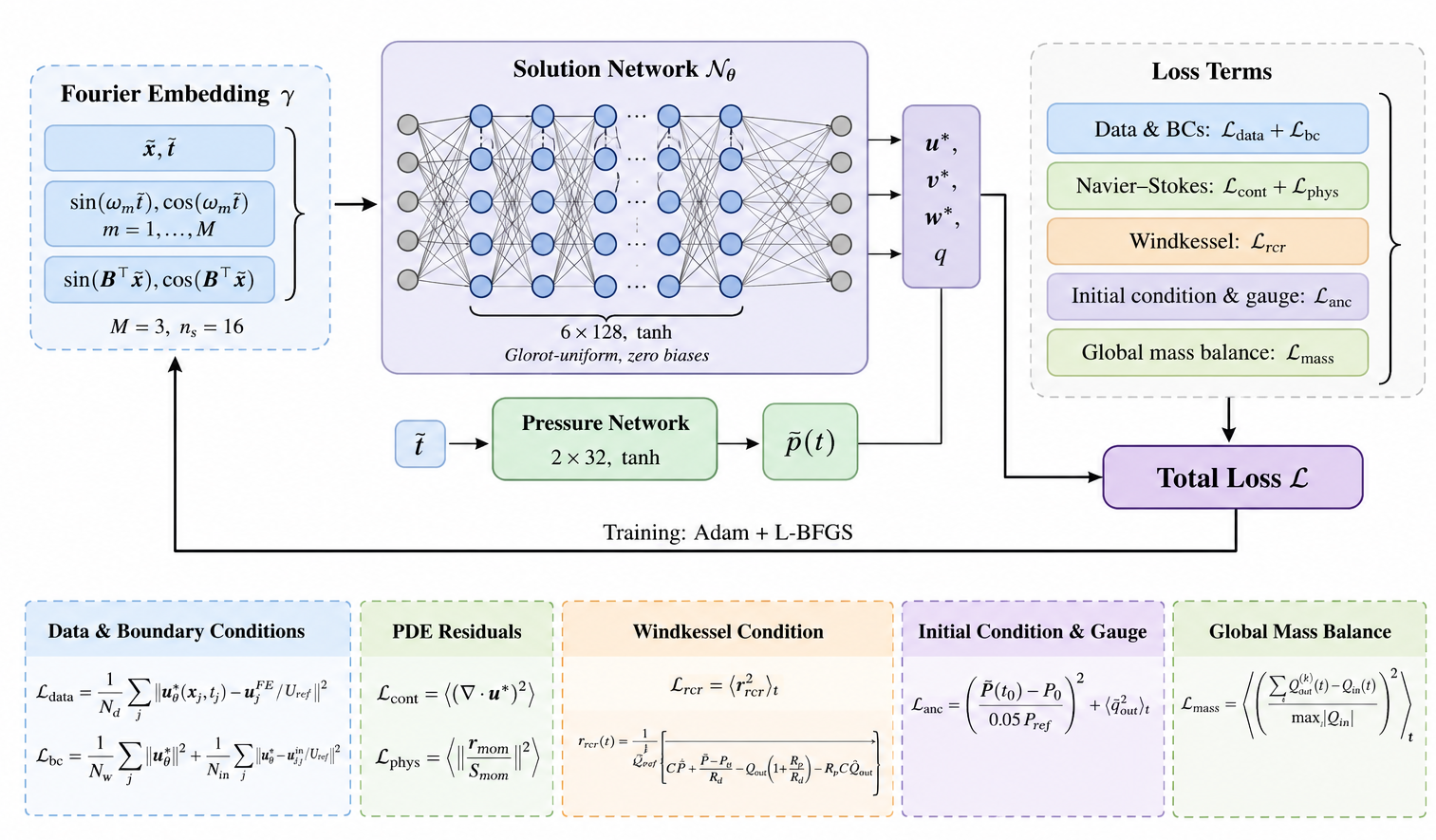}
    \caption{Schematic PINN architecture, including the Fourier embedding, the solutions and pressure network, and the main components of the composite loss function.}
    \label{fig:pinn_arch}
\end{figure*}

This section describes the PINN architecture and approach used throughout the paper, illustrated schematically in Fig.~\ref{fig:pinn_arch}. The network approximates the velocity and pressure fields over the space-time domain $\Omega\times[0,T]$ and is trained against the Navier--Stokes equations of Section \ref{subsec:navier_stokes}, the no-slip and inflow conditions, the RCR relation of Eq. \eqref{eq:rcr_ode} at every outlet, and reference data generated with finite element simulations (FEBio, fluid module; \cite{maas2012febio}). The implementation is written in PyTorch \cite{paszke2019pytorch}.
 
\subsection{Problem set-up and training data}\label{subsec:pinn_data}
 
The reference simulation provides, at the mesh nodes $\mathbf{x}_i$ and at a set of output times $t_k$, the nodal velocity $\mathbf{u}^{\mathrm{FE}}_{ik}$ and the element pressure, which is averaged onto the nodes to obtain $p^{\mathrm{FE}}_{ik}$. Four families of points are built from the mesh:
\begin{itemize}
    \item \emph{Interior collocation points}, drawn at random inside every element (barycentric weights for tetrahedra, convex weights for hexahedra) and at random times, where the Navier--Stokes residuals are evaluated.
    \item \emph{Wall and inlet points}, taken from the mesh surfaces. At the wall $\mathbf{u}=\mathbf{0}$; at the inlet the FE nodal velocity is interpolated in time with a monotone cubic (PCHIP) interpolant, so that the boundary values are available at arbitrary times and not only at the (non-uniform) output times of the solver.
    \item \emph{Data points}, a random fraction of the FE nodes at the output times, where the velocity (and, optionally, the pressure) is supervised.
    \item \emph{Outlet points}, all the nodes of every outlet surface, at a uniform set of $N_t^{\mathrm{rcr}}$ instants covering $[0,T]$, used to evaluate the surface integrals of the Windkessel model.
\end{itemize}
A disjoint set of FE nodes, never used for training, is reserved for validation.
 
Surface integrals over an outlet $\Gamma_{out}^{(k)}$ are discretized with lumped nodal areas $A_i$ (one quarter of the area of each adjacent quadrilateral, one third for triangles) and area-weighted vertex normals $\mathbf{n}_i$ computed directly from the surface mesh and oriented outwards. Hence
\begin{equation}\label{eq:Qdisc}
    Q_{out}^{(k)}(t)\approx\sum_{i\in\Gamma_{out}^{(k)}} A_i\,\mathbf{u}(\mathbf{x}_i,t)\cdot\mathbf{n}_i, \qquad
    \bar P^{(k)}(t)\approx\frac{\sum_i A_i\,p(\mathbf{x}_i,t)}{\sum_i A_i}.
\end{equation}
The corresponding inflow $Q_{in}(t)$ is obtained in the same way from the inlet surface with inward orientation.
 
\subsection{Non-dimensionalization}\label{subsec:pinn_scaling}
 
Spatial coordinates and time are first mapped to $[-1,1]$ with the bounding box of the domain,
\begin{equation}
    \tilde x_j=\frac{2\,(x_j-x_j^{\min})}{x_j^{\max}-x_j^{\min}}-1,\qquad
    \tilde t=\frac{2\,(t-t^{\min})}{t^{\max}-t^{\min}}-1 ,
\end{equation}
and these normalized variables are the network inputs. Any axis whose extent is smaller than $5\%$ of the longest one (as in the thin slab used for the canonical test) is widened to the longest extent before normalizing; otherwise the chain-rule factor $2L_{ref}/(x_j^{\max}-x_j^{\min})$ becomes very large along that axis and makes the residual extremely stiff.
 
The physical variables are scaled with a reference velocity $U_{ref}$ (maximum inlet speed), a reference length $L_{ref}$ (vessel diameter or longest dimension of the domain), $T_{ref}=L_{ref}/U_{ref}$, and two pressure scales. The pressure is decomposed as
\begin{equation}\label{eq:pdecomp}
    p(\mathbf{x},t)=P_{ref}\,\bar p(t)+P_{dyn}\,q(\mathbf{x},t),
\end{equation}
where $\bar p(t)$ is a spatially uniform basal pressure level, scaled with $P_{ref}=\max|p^{\mathrm{FE}}|$, and $q$ is the spatially varying part, scaled with $P_{dyn}=\max|p^{\mathrm{FE}}-\langle p^{\mathrm{FE}}\rangle_\Omega|$. This separation is the key to the pressure-scale disparity discussed in the Introduction: the momentum equation only involves $\nabla p=P_{dyn}\nabla q$, so the large basal pressure imposed by the Windkessel never competes with the small local gradients that drive the flow; both are represented by outputs of order one. Denoting by an asterisk the non-dimensional quantities, the governing equations read
\begin{align}
    \nabla^*\!\cdot\mathbf{u}^* &= 0,\label{eq:nd_cont}\\
    \frac{\partial\mathbf{u}^*}{\partial t^*}+(\mathbf{u}^*\!\cdot\nabla^*)\mathbf{u}^*+\kappa_p\nabla^* q &=\frac{1}{Re}\nabla^{*2}\mathbf{u}^*,\label{eq:nd_mom}
\end{align}
with $Re=\rho U_{ref}L_{ref}/\mu$ and $\kappa_p=P_{dyn}/(\rho U_{ref}^2)$. In viscous-dominated problems with small pressure variations all the terms of Eq. \eqref{eq:nd_mom} are much smaller than one, so the momentum residual is divided by a characteristic magnitude $S_{mom}=\kappa_p\,L_{ref}/(0.5\,L_{max})$, with $L_{max}$ the longest extent of the domain; this keeps the physics loss comparable to the data loss.
 
\subsection{Network architecture}\label{subsec:pinn_arch}
 
The solution fields are represented by a fully connected network $\mathcal{N}_\theta$ with $6$ hidden layers of $128$ neurons and hyperbolic tangent activation, Glorot-uniform initialization and zero biases,
\begin{equation}
    \big(u^*,v^*,w^*,q\big)=\mathcal{N}_\theta\!\big(\boldsymbol{\gamma}(\tilde{\mathbf{x}},\tilde t)\big),
\end{equation}
preceded by a fixed (non-trainable) Fourier embedding $\boldsymbol{\gamma}$ that is evaluated inside the forward pass, so that automatic differentiation with respect to the raw inputs remains exact:
\begin{equation}
    \boldsymbol{\gamma}=\Big[\tilde{\mathbf{x}},\;\tilde t,\;\sin(\omega_m\tilde t),\cos(\omega_m\tilde t)\big|_{m=1}^{M},\;\sin(\mathbf{B}^{\!\top}\tilde{\mathbf{x}}),\cos(\mathbf{B}^{\!\top}\tilde{\mathbf{x}})\Big].
\end{equation}
The temporal frequencies are the first $M$ harmonics of the cardiac period $T_c$, $\omega_m=\pi m\,(t^{\max}-t^{\min})/T_c$ (rad per unit of $\tilde t$), where $T_c$ is detected automatically from the dominant peak of the zero-padded FFT of the inflow $Q_{in}(t)$. This embedding counteracts the spectral bias of $\tanh$ networks \cite{rahaman2019spectral, tancik2020fourier}, which, with raw time as the only input, tend to return a single smooth hump instead of several cardiac cycles. The matrix $\mathbf{B}\in\mathbb{R}^{3\times n_s}$ contains random low-frequency Gaussian features ($n_s=16$, standard deviation $1$) that enrich the spatial representation without introducing high-frequency noise. We use $M=3$.
 
The basal pressure $\bar p(t)$ of Eq. \eqref{eq:pdecomp} is produced by a second, small network (two hidden layers of $32$ neurons) that receives $\tilde t$ and its periodic features. Because it depends only on time, the non-periodic character of the pressure (e.g., the charging of the compliance from rest) can be represented separately from the periodic velocity dynamics.
 
\subsection{Composite loss function}\label{subsec:pinn_loss}
 
The network is trained by minimizing a weighted sum of eight terms. All of them are mean squared errors of non-dimensional quantities.
 
\paragraph{Data and boundary conditions.} With $N_d$ data points, $N_w$ wall points and $N_{in}$ inlet points,
\begin{align}
    \mathcal{L}_{data}&=\frac{1}{N_d}\sum_{j}\big\|\mathbf{u}^*_\theta(\mathbf{x}_j,t_j)-\mathbf{u}^{\mathrm{FE}}_j/U_{ref}\big\|^2,\\
    \mathcal{L}_{bc}&=\frac{1}{N_w}\sum_j\|\mathbf{u}^*_\theta\|^2+\frac{1}{N_{in}}\sum_j\big\|\mathbf{u}^*_\theta-\mathbf{u}^{\mathrm{in}}_j/U_{ref}\big\|^2 .
\end{align}
An optional pressure data term, $\mathcal{L}_{p}=\langle((p_\theta-p^{\mathrm{FE}})/P_{dyn})^2\rangle$, is available; in the forward problems of this paper its weight is zero, so the pressure is \emph{inferred} from the physics and the outlet model and not supervised.
 
\paragraph{Navier--Stokes residuals.} Evaluated at the interior collocation points with first and second derivatives obtained by automatic differentiation and rescaled to non-dimensional variables by the factors $2L_{ref}/(x_j^{\max}-x_j^{\min})$ and \\ $2T_{ref}/(t^{\max}-t^{\min})$,
\begin{equation}
    \mathcal{L}_{cont}=\big\langle(\nabla^*\!\cdot\mathbf{u}^*)^2\big\rangle,\qquad
    \mathcal{L}_{phys}=\Big\langle\big\|\mathbf{r}_{mom}/S_{mom}\big\|^2\Big\rangle,
\end{equation}
where $\mathbf{r}_{mom}$ is the left-minus-right side of Eq. \eqref{eq:nd_mom}.
 
\paragraph{Windkessel outlet condition.} For each outlet $k$, $Q_{out}^{(k)}(t)$ and $\bar P^{(k)}(t)$ are computed from Eq. \eqref{eq:Qdisc} at the $N_t^{\mathrm{rcr}}$ instants. Since points at different times are independent, the time derivative of a surface integral is the sum over nodes of the pointwise time derivative, so $\dot Q_{out}$ and $\dot{\bar P}$ require only one backward pass each. The residual of Eq. \eqref{eq:rcr_ode}, including a possible offset pressure $P_d$, is
\begin{equation}\label{eq:rcr_res}
    r_{rcr}(t)=\frac{1}{Q_{ref}}\Big[C\dot{\bar P}+\frac{\bar P-P_d}{R_d}-Q_{out}\Big(1+\frac{R_p}{R_d}\Big)-R_pC\,\dot Q_{out}\Big],
\end{equation}
with $Q_{ref}=U_{ref}A_{out}$, and $\mathcal{L}_{rcr}=\langle r_{rcr}^2\rangle_t$ averaged over time and outlets. The parameters $R_p$, $R_d$, $C$ and $P_d$ are read from the FE model, so that the PINN and the reference solution share exactly the same outlet condition.
 
\paragraph{Initial condition and gauge.} Equation \eqref{eq:rcr_ode} is a first-order ODE, so its residual alone does not determine the pressure level. We add
\begin{equation}
    \mathcal{L}_{anc}=\Big(\frac{\bar P(t_0)-P_0}{0.05\,P_{ref}}\Big)^2+\big\langle\bar q_{out}(t)^2\big\rangle_t ,
\end{equation}
where $P_0$ is the initial outlet pressure and $\bar q_{out}$ is the area-averaged $q$ over the outlet. The second term fixes the gauge of the decomposition: the area-averaged outlet pressure equals $P_{ref}\bar p(t)$ exactly, so $\bar p(t)$ is the outlet pressure seen by the Windkessel.
 
\paragraph{Global mass balance.} The local residual $\langle(\nabla\!\cdot\mathbf{u})^2\rangle$ is small even when a sizeable fraction of the inflow fails to leave the domain. We therefore impose the integral constraint
\begin{equation}
    \mathcal{L}_{mass}=\Big\langle\Big(\frac{\sum_k Q_{out}^{(k)}(t)-Q_{in}(t)}{\max_t|Q_{in}|}\Big)^2\Big\rangle_t ,
\end{equation}
where $Q_{in}(t)$ is computed from the inlet velocity data and interpolated at the instants of the outlet set.
 
\paragraph{Total loss.} The terms are combined as
\begin{equation}\label{eq:total_loss}
\begin{aligned}
    \mathcal{L} &= \lambda_{d}\mathcal{L}_{data} + \lambda_{b}\mathcal{L}_{bc} + \lambda_{m}\mathcal{L}_{mass} + \lambda_{a}\mathcal{L}_{anc} \\
    &\quad + r(e)\big[\lambda_{p}\mathcal{L}_{phys} + \lambda_{c}\mathcal{L}_{cont} + \lambda_{r}\mathcal{L}_{rcr}\big].
\end{aligned}
\end{equation}
where the ramp $r(e)=\min(1,e/E_r)$ with $E_r=1500$ epochs switches the PDE and Windkessel residuals on gradually, so that the network first learns the data and the boundary values and the stiff terms are introduced on a reasonable starting point. The weights are listed in Table \ref{tbl:hyper}.
 
\subsection{Optimization strategy}\label{subsec:pinn_opt}
 
Training proceeds in two stages. In the first stage the Adam optimizer is run for $E=8000$ epochs with a linear warm-up over $200$ epochs followed by an exponential decay of the learning rate from $10^{-3}$ to $10^{-5}$, and with gradient-norm clipping at $5$. At each epoch a random mini-batch is drawn from every point set (Table \ref{tbl:hyper}), whereas the Windkessel, mass and anchor terms use all outlet nodes at all $N_t^{\mathrm{rcr}}$ instants. In the second stage the solution is polished with up to three rounds of L-BFGS \cite{liu1989limited} (strong-Wolfe line search, history of $50$), each with a new fixed collocation set of twice the mini-batch size. After every round the relative $L_2$ velocity error on the held-out nodes is evaluated; a round that does not improve it is discarded and the previous parameters are restored. During training we monitor the relative $L_2$ error of the velocity, the RMSE of the dynamic pressure $(p-\langle p\rangle_\Omega)$ and of the absolute pressure on the held-out nodes, and the relative RMSE between $\sum_kQ_{out}^{(k)}$ and $Q_{in}$.

\begin{table}[t]
\caption{Hyper-parameters of the PINN.}\label{tbl:hyper}
\small
\begin{tabular*}{\tblwidth}{@{}LL@{}}
\toprule
Parameter & Value \\
\midrule
Network & $6\times128$, $\tanh$; pressure net $2\times32$ \\
 Embedding & $M=3$ harmonics of $T_c$; \\
& $n_s=16$ spatial features \\
Loss weights & $\lambda_d=10$, $\lambda_b=10$, $\lambda_m=5$, $\lambda_a=1$, \\
             & $\lambda_p=\lambda_c=\lambda_r=1$ \\
Ramp & $E_r=1500$ epochs \\
Adam & $E=8000$; $\eta: 10^{-3} \to 10^{-5}$; clip $5$ \\
L-BFGS & $3$ rounds $\times$ $200$ iterations \\
Mini-batch & $20000$ interior, $4096$ boundary, $8192$ data \\
Windkessel times & $N_t^{\mathrm{rcr}}=120$ \\
Data & $15\%$ of FE nodes, all output times; \\
     & $400$ validation nodes \\
\bottomrule
\end{tabular*}
\end{table}

\section{Validation on Canonical Two-Dimensional Geometries}\label{sec:validation}

Prior to addressing the patient-specific 3D geometry, the proposed PINN framework is systematically verified on three canonical 2D configurations of increasing physical and geometrical complexity: (i) a straight channel under steady laminar flow with an analytical plane Poiseuille solution; (ii) a channel with a lateral offset (double bend), which introduces flow acceleration, non-uniform pressure gradients, and local geometric discontinuities; and (iii) a straight channel under transient pulsatile flow coupled to an outlet 3-element Windkessel (RCR) model. In all cases, reference numerical solutions are generated using FEBio on quasi-2D meshes, and the PINN predictions are evaluated on held-out validation nodes not used during network training.

\subsection{Error metrics}\label{subsec:metrics}

The agreement between the PINN predictions and the FEBio reference solution is quantified using the relative $L_2$ velocity error ($\varepsilon_{L_2}$), the root-mean-square error ($\mathrm{RMSE}$), and the coefficient of determination ($R^2$), evaluated over the set of held-out validation nodes $\mathcal{V}$:
\begin{equation}\label{eq:metrics_l2_rmse}
\begin{aligned}
    \varepsilon_{L_2} &= \frac{\big(\sum_{\mathcal{V}}\|\mathbf{u}_\theta-\mathbf{u}^{\mathrm{FE}}\|^2\big)^{1/2}}{\big(\sum_{\mathcal{V}}\|\mathbf{u}^{\mathrm{FE}}\|^2\big)^{1/2}}, \\[1.5ex]
    \mathrm{RMSE} &= \left(\frac{1}{|\mathcal{V}|}\sum_{\mathcal{V}}\|\mathbf{u}_\theta-\mathbf{u}^{\mathrm{FE}}\|^2\right)^{1/2}.
\end{aligned}
\end{equation}
\begin{equation}
    R^2=1-\frac{\sum_{\mathcal{V}}\|\mathbf{u}_\theta-\mathbf{u}^{\mathrm{FE}}\|^2}{\sum_{\mathcal{V}}\|\mathbf{u}_\theta-\bar{\mathbf{u}}^{\mathrm{FE}}\|^2}.
\end{equation}
Because the incompressible Navier--Stokes formulation only involves the pressure gradient $\nabla p$, the predicted pressure is assessed through its dynamic component, $p-\langle p\rangle_\Omega$, as well as its absolute level. Global mass conservation is evaluated via the normalized flow rate ratio $Q(x)/Q_{in}$ across transversal cross-sections along the domain length.

\subsection{Straight channel under steady flow}\label{subsec:straight_steady}

The first benchmark considers a straight channel of length $L=10$~mm and height $h=3$~mm with a parabolic velocity profile of mean magnitude $U$ prescribed at the inlet and zero reference pressure at the outlet. For fully developed laminar flow, the velocity profile and spatial pressure gradient follow the analytical plane Poiseuille solution:
\begin{equation}
    u(y)=\frac{3}{2}\,U\Big(1-\frac{4y^2}{h^2}\Big),\qquad
    \frac{\mathrm{d}p}{\mathrm{d}x}=-\frac{12\,\mu\,U}{h^2},
\end{equation}
yielding a theoretical total pressure drop of $\Delta P = 12\mu U L/h^2$. Two physical regimes were evaluated: $\mu=0.004$~Pa$\cdot$s with $U=0.7$~m/s (analytical drop of $37.3$~Pa) and $\mu=0.005$~Pa$\cdot$s with $U=0.3$~m/s ($20.0$~Pa).  Two network variants were trained: one supervised exclusively with sparse velocity data, and one incorporating explicit pressure supervision (hereafter referred to as \emph{with $p$} and \emph{without $p$}).

\begin{figure*}[htbp]
    \centering

    \begin{minipage}[t]{0.48\textwidth}
        \centering
        \includegraphics[width=1.05\linewidth]{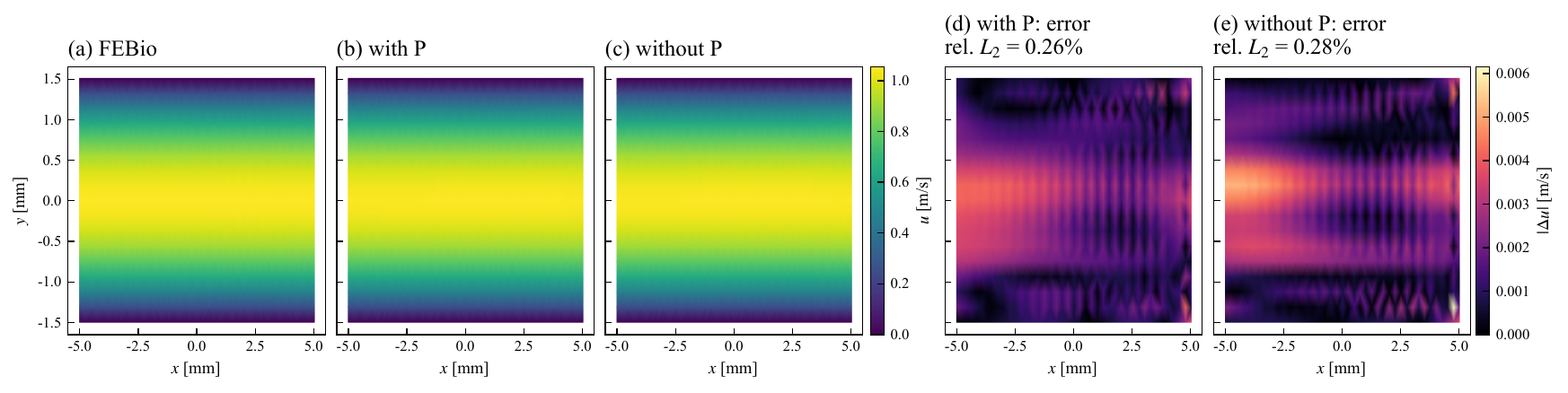}

    \end{minipage}
    \hfill
    \begin{minipage}[t]{0.48\textwidth}
        \centering
        \includegraphics[width=1.05\linewidth]{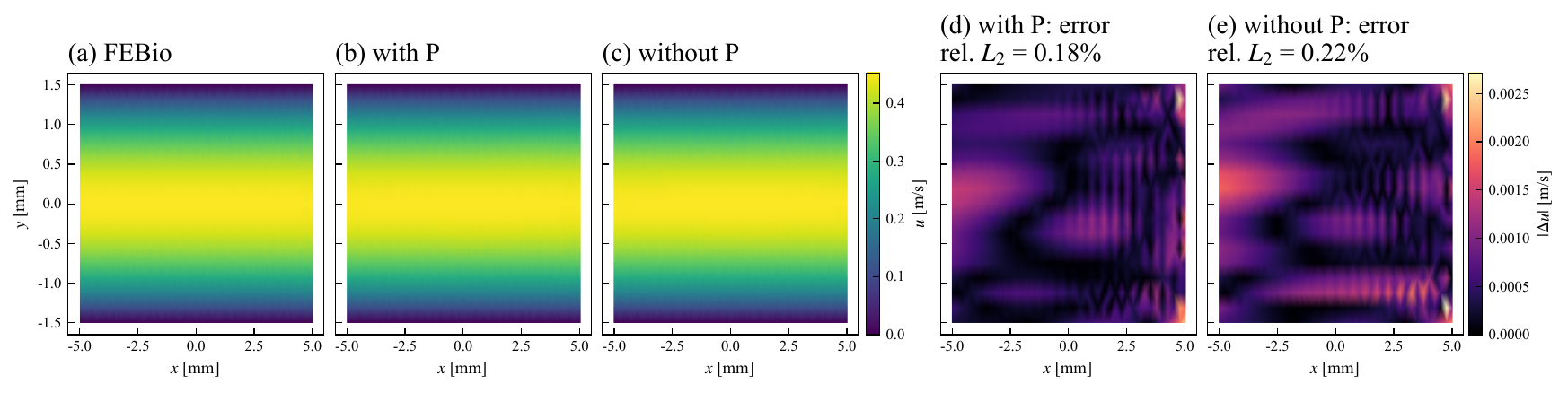}
        
    \end{minipage}

    \vspace{0.3cm}

    \begin{minipage}[t]{0.48\textwidth}
        \centering
        \includegraphics[width=1.05\linewidth]{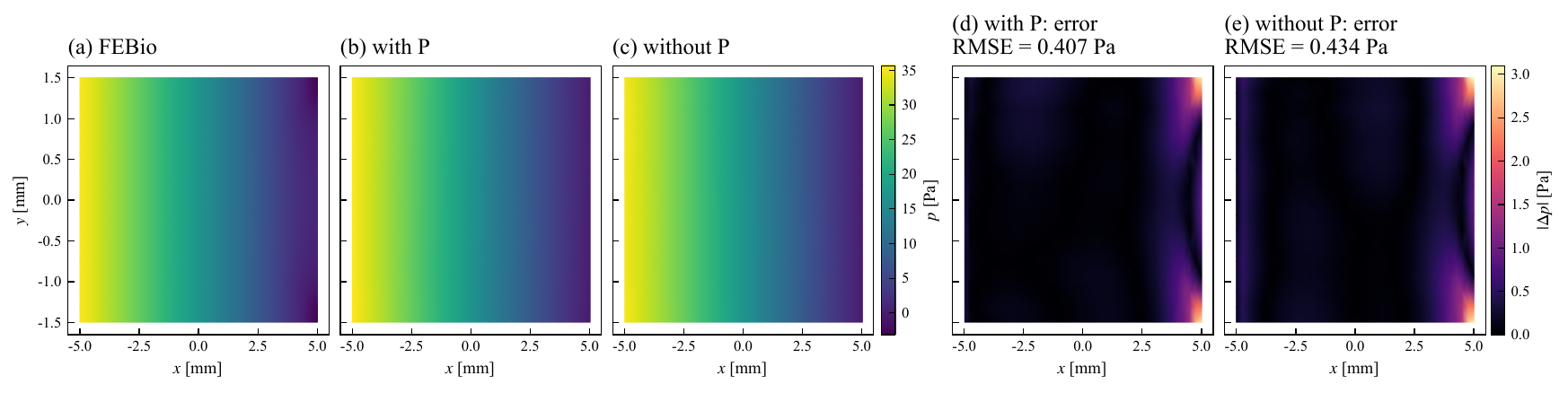}

    \end{minipage}
    \hfill
    \begin{minipage}[t]{0.48\textwidth}
        \centering
        \includegraphics[width=1.05\linewidth]{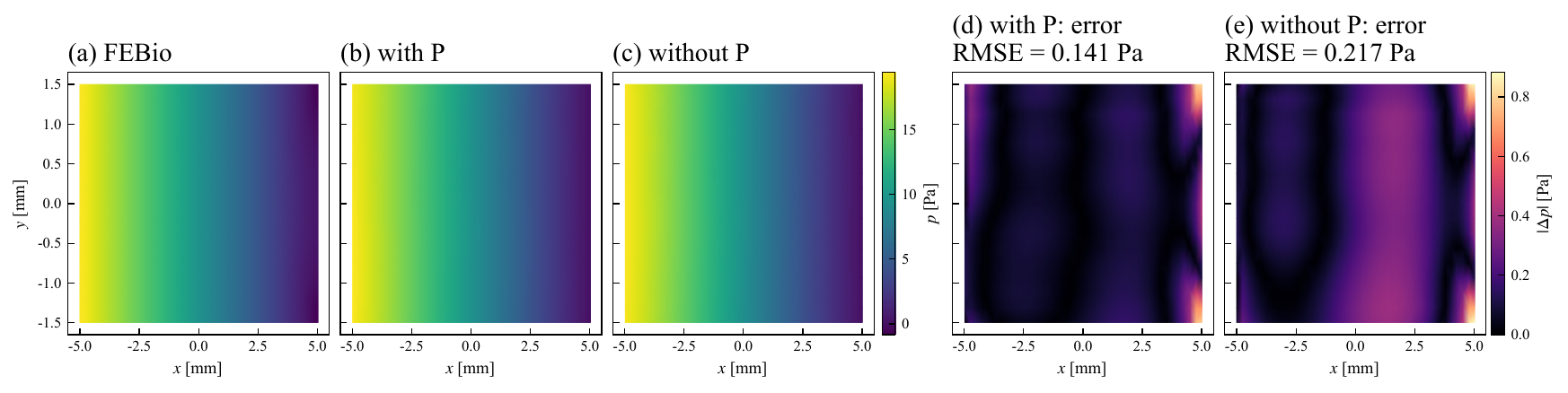}
        
    \end{minipage}

    \vspace{0.3cm}
    \begin{minipage}[t]{0.48\textwidth}
        \centering
        \includegraphics[width=\linewidth]{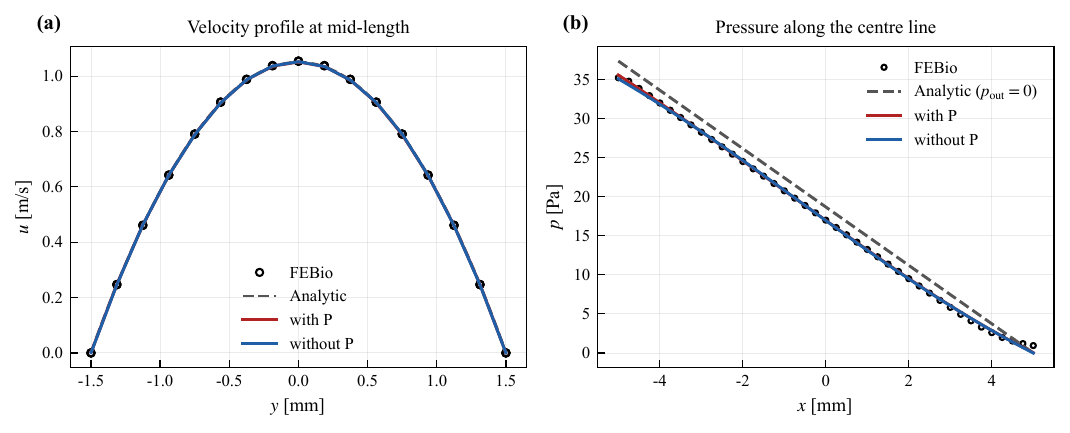}

    \end{minipage}
    \hfill
    \begin{minipage}[t]{0.48\textwidth}
        \centering
        \includegraphics[width=\linewidth]{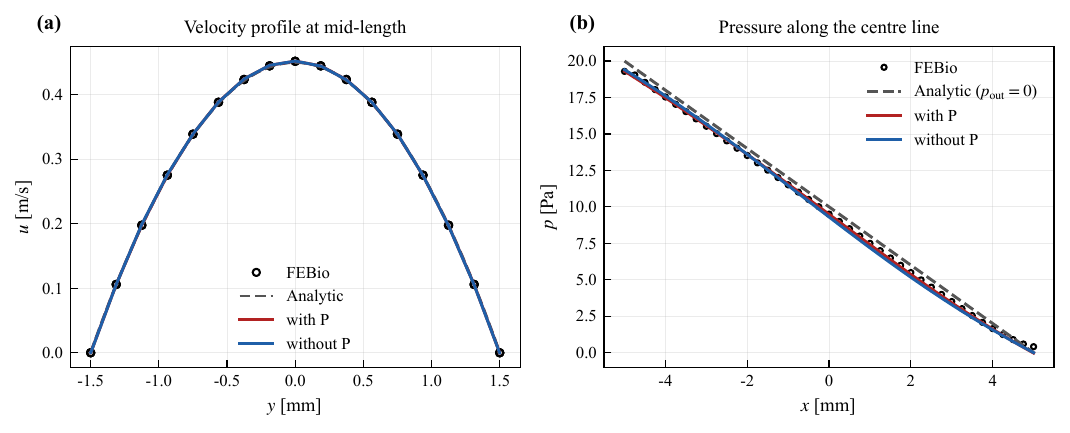}

    \end{minipage}

    \caption{Validation, plane Poiseuille flow. Left: Case 1 ($\mu$ = 0.004 Pa s, $U_{in}$ = 0.7 m/s). Right: Case 2 ($\mu$ = 0.005 Pa s, $U_{in}$ = 0.3 m/s). (a) Velocity profile at mid-length and (b) pressure along the centre line: FEBio (circles), analytic solution (dashed) and PINN (solid). Fields: FEBio and PINN share the same colour scale; the right/bottom maps show the absolute PINN-FEBio error.
}
    \label{fig:straight_steady}
\end{figure*}

Figure~\ref{fig:straight_steady} demonstrates that the predicted mid-length velocity profiles match both the FEBio reference and the analytical parabola, reaching peak centerline velocities of $1.05$~m/s and $0.45$~m/s, respectively. The maximum pointwise velocity error remains below $0.6\%$ of the peak speed ($4.1\times10^{-3}$~m/s for Case 1 and $2.5\times10^{-3}$~m/s for Case 2). The linear pressure drop along the centerline is accurately captured. The slight deviation between the numerical solvers and the theoretical line ($\approx35$~Pa vs. $37.3$~Pa) is shared by FEBio and is attributed to finite entrance length effects. Pressure errors remain negligible in the domain interior, peaking slightly near the outlet boundary where $p=0$ is imposed. Furthermore, both network variants yield nearly identical predictions, confirming that pressure supervision is redundant in simple canonical geometries.

\begin{table}[htbp]
\centering
\caption{Training and validation cases for the straight-channel benchmark.}
\label{tab:straight_steady_cases}
\small
\begin{tabular}{c c c c}
\hline
$\mu$ & Training $U_{\mathrm{in}}$ & Validation $U_{\mathrm{in}}$ & $N$ \\
(Pa$\cdot$s) & (m/s) & (m/s) & \\
\hline
0.003 & 0.1--0.9 & -- & 5 \\
0.004 & 0.1, 0.3, 0.5, 0.9 & 0.7 & 5 \\
0.005 & 0.1, 0.5, 0.7, 0.9 & 0.3 & 5 \\
\hline
Total & 13 & 2 & 15 \\
\hline
\end{tabular}
\end{table}

The dataset comprises 15 parameter combinations, of which 13 were used for training and two were reserved for validation. The validation cases correspond to $(\mu,U_{\mathrm{in}})=(0.004,0.7)$~Pa$\cdot$s and m/s, respectively, and $(0.005,0.3)$~Pa$\cdot$s and m/s, respectively. All four corner combinations of the parameter space were included in the training set. Therefore, the validation assesses interpolation within the sampled parameter range rather than extrapolation beyond it.

\FloatBarrier
\subsection{Channel with lateral offset under steady flow}\label{subsec:curved}

Geometric complexity is introduced using a channel of width $20$~mm and total length $130$~mm, whose longitudinal axis is shifted laterally by $40$~mm between $x\approx45$~mm and $x\approx85$~mm. The double bend generates local flow acceleration, internal streamline curvature, non-uniform pressure fields, and localized velocity gradients.

A total of 16 parameter combinations were considered, spanning four viscosity values ($\mu=0.002$, $0.003$, $0.004$, and $0.005$~Pa$\cdot$s) and four inlet velocities ($U_{\mathrm{in}}=0.05$, $0.1$, $0.15$, and $0.2$~m/s). Fourteen cases were used for training, while two were reserved for validation (Table~\ref{tab:curved_steady_cases}). The validation cases were $(\mu,U_{\mathrm{in}})=(0.003,0.15)$ and $(0.004,0.1)$, respectively. Both correspond to combinations within the training parameter range and therefore assess interpolation performance.

\begin{table}[htbp]
\centering
\caption{Training and validation cases for the channel with lateral offset.}
\label{tab:curved_steady_cases}
\small
\begin{tabular}{c c c c}
\hline
$\mu$ & Training $U_{\mathrm{in}}$ & Validation $U_{\mathrm{in}}$ & $N$ \\
(Pa$\cdot$s) & (m/s) & (m/s) & \\
\hline
0.002 & 0.05, 0.1, 0.15, 0.2 & -- & 4 \\
0.003 & 0.05, 0.1, 0.2 & 0.15 & 4 \\
0.004 & 0.05, 0.15, 0.2 & 0.1 & 4 \\
0.005 & 0.05, 0.1, 0.15, 0.2 & -- & 4 \\
\hline
Total & 14 & 2 & 16 \\
\hline
\end{tabular}
\end{table}

The objective of this benchmark is to assess whether the trained PINN can reproduce the flow fields obtained with FEBio for previously unseen parameter combinations, while preserving the expected physical behaviour of the flow. In particular, the comparison focuses on velocity and pressure distributions, as well as global mass conservation.

\begin{figure*}[htbp]
    \centering

    \begin{minipage}[t]{0.48\textwidth}
        \centering
        \includegraphics[width=\linewidth]{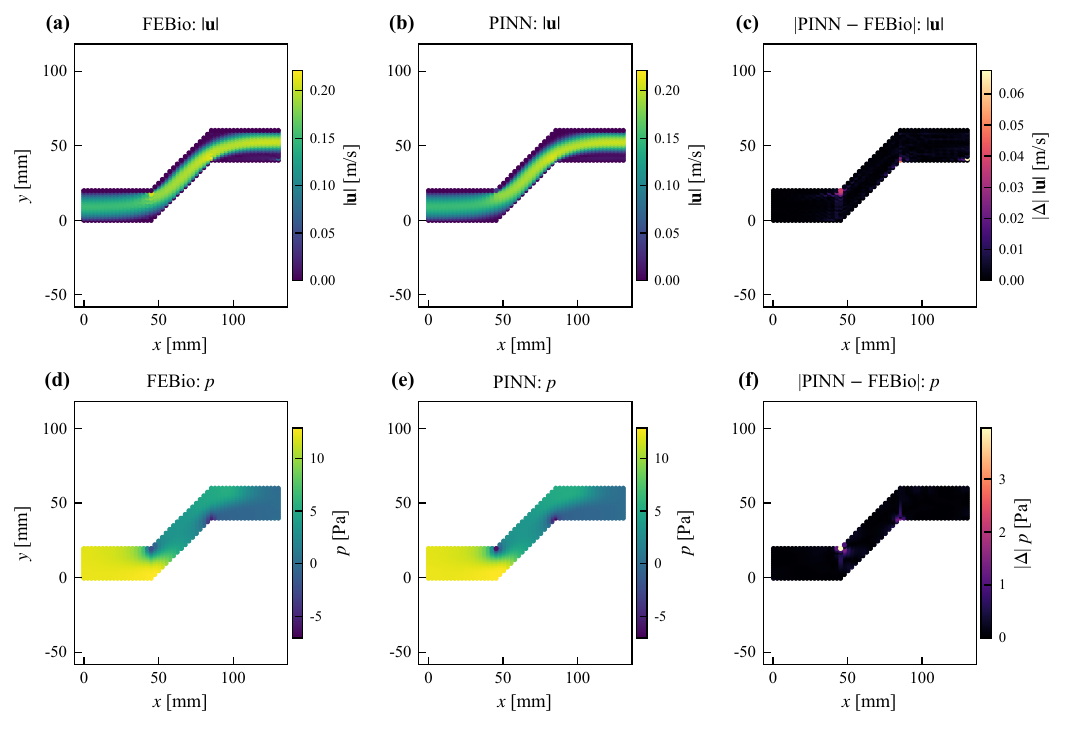}

    \end{minipage}
    \hfill
    \begin{minipage}[t]{0.48\textwidth}
        \centering
        \includegraphics[width=\linewidth]{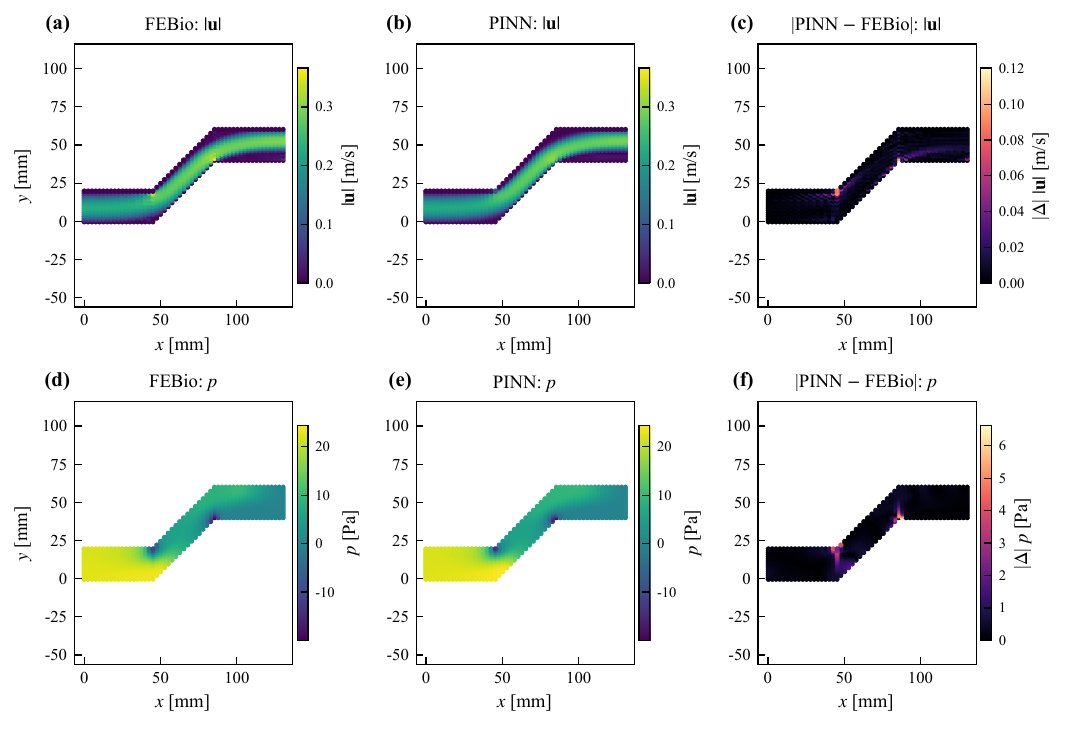}

    \end{minipage}

    \vspace{0.3cm}

    \begin{minipage}[t]{0.48\textwidth}
        \centering
        \includegraphics[width=\linewidth]{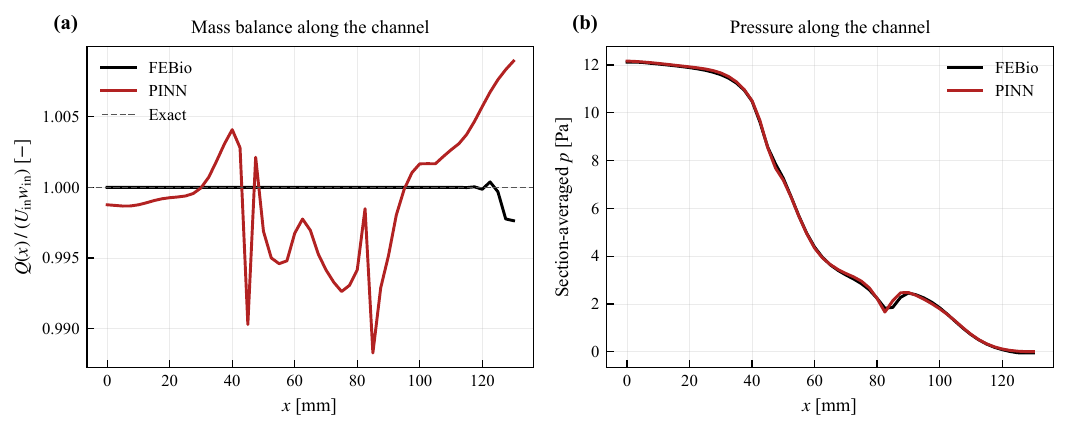}

    \end{minipage}
    \hfill
    \begin{minipage}[t]{0.48\textwidth}
        \centering
        \includegraphics[width=\linewidth]{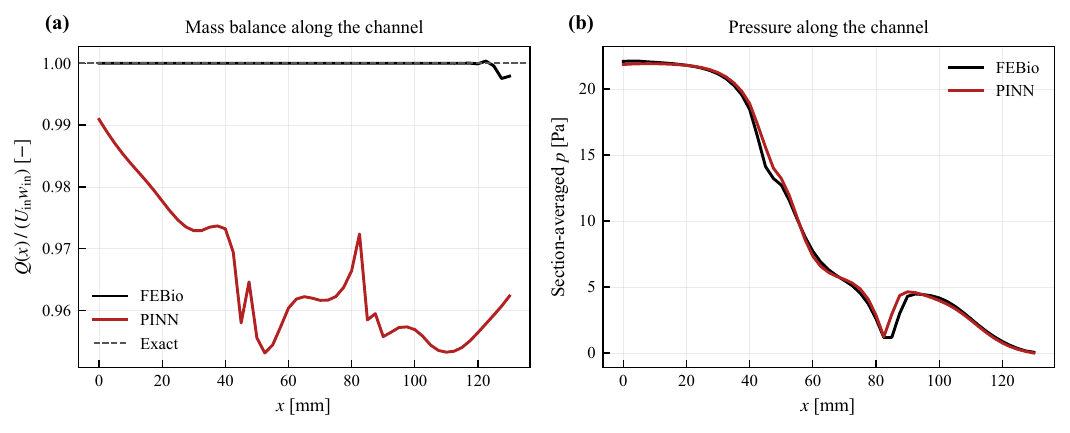}

    \end{minipage}

    \caption{Validation, Z-bend channel. Left: Case 1 ($\mu$ = 0.004 Pa s, $U_{in}$ = 0.1 m/s).  Right: Case 2 ($\mu$ = 0.003 Pa s, $U_{in}$ = 0.15 m/s). Fields: FEBio and PINN share the same colour scale per variable; the right column shows the absolute PINN-FEBio error. (a-f) velocity magnitude and pressure. Lower figure: flux $Q(x)$ normalised by the inlet flux (exact value 1) and section-averaged pressure for FEBio and PINN.
}
    \label{fig:z_bend_channel}
\end{figure*}

Figure~\ref{fig:z_bend_channel} presents the results for the two validation cases: $\mu=0.003$~Pa$\cdot$s with $U_{\mathrm{in}}=0.15$~m/s (peak velocity $0.35$~m/s and inlet pressure approximately $22$~Pa), and $\mu=0.004$~Pa$\cdot$s with $U_{\mathrm{in}}=0.1$~m/s (peak velocity $0.20$~m/s and inlet pressure approximately $12$~Pa).

Results are summarized in Fig.~\ref{fig:z_bend_channel}. Across both cases, the PINN accurately reconstructs the full velocity magnitude field, capturing local acceleration through the inclined neck and downstream flow re-development. Velocity errors remain negligible along straight segments, concentrating at sharp inner corners ($0.12$~m/s and $0.06$~m/s peak pointwise errors). Pressure errors exhibit a similar pattern, with localized peaks at corners ($\approx30\%$ local error) and negligible values throughout the rest of the channel.

Section-averaged pressure profiles along the channel axis highlight the role of pressure data. When sparse pressure supervision is provided, the PINN tracks the FEBio reference across the entire domain, including the localized pressure minimum at $x\approx82$~mm. Without pressure supervision, the spatial gradient trend is preserved, but an absolute pressure level shift occurs ($\approx0.7$--$0.9$~Pa offset). Global mass conservation $Q(x)/Q_{in}$ is satisfied by FEBio within $0.3\%$, whereas the soft continuity penalty in the PINN allows a maximum flow deficit of $4.7\%$ ($6.7\%$ without $p$) for the higher-velocity case. These observations justified the inclusion of the explicit global mass balance loss term $\mathcal{L}_{mass}$ in the final architecture.

Global mass conservation is satisfied by FEBio within $0.3\%$, whereas the PINN exhibits a maximum flow deficit of $4.7\%$ for the higher-velocity case, increasing to $6.7\%$ when pressure supervision is excluded. Although the PINN does not enforce perfectly constant flow rate along the channel, the predicted flux remains reasonably close to the inlet value. This residual variation reflects the limitations of the soft continuity constraint, particularly in geometries with strong local flow gradients. These results motivated the inclusion of an explicit global mass balance loss term, $\mathcal{L}_{\mathrm{mass}}$, in the final architecture.

 \FloatBarrier
\section{Straight channel under pulsatile flow with RCR outlet}\label{subsec:straight_pulsatile}

To evaluate transient performance and dynamic lumped-parameter coupling, the straight channel ($10\times3$~mm) is subjected to a periodic inflow profile ($T=1$~s period over two cardiac cycles) and coupled to a 3-element Windkessel (RCR) model at the outlet. The peak centerline velocity reaches $\approx1.5$~m/s.

The FE model is a quasi-2D slab of $10\times3\times0.1$~mm discretized with $3200$ hexahedral elements and $4182$ nodes ($41\times17\times6$), with a parabolic inlet velocity whose time dependence follows $\tfrac12(1-\cos2\pi t)$, no-slip walls, $\rho=1060$~kg/m$^3$ and $\mu=0.004$~Pa$\cdot$s. The resulting RCR time constant is $R_dC\approx0.49$~s.

\begin{table}[!ht]
\caption{Parameters of the 3-element Windkessel (RCR) model coupled at the outlet of the 2D pulsatile straight channel benchmark (read from the FEBio input file).}\label{tbl:rcr_2d}
\begin{tabular*}{\tblwidth}{@{}LL@{}}
\toprule
Parameter & Value \\
\midrule
Proximal resistance ($R_p$) & $3.5\times10^{7}$ $\mathrm{Pa\cdot s/m^3}$ \\
Distal resistance ($R_d$) & $4.62\times10^{7}$ $\mathrm{Pa\cdot s/m^3}$ \\
Arterial compliance ($C$) & $1.07\times10^{-8}$ $\mathrm{m^3/Pa}$ \\
Offset pressure ($P_d$) & $0$ $\mathrm{Pa}$ \\
Initial pressure & $0$ $\mathrm{Pa}$ \\
\bottomrule
\end{tabular*}
\end{table}

To assess generalisation, the pulsatile straight channel of Section~\ref{subsec:straight_pulsatile} is extended to a
family of FEBio simulations that share the geometry, the mesh, the period ($T=1$~s, two cycles) and the RCR outlet
(Table~\ref{tbl:rcr_2d}) and differ in the dynamic viscosity $\mu$ and in the peak mean inlet velocity $U$ (the inflow is
$U\,\tfrac12(1-\cos 2\pi t)$ with a parabolic profile). A single PINN with inputs $(x,y,t,\mu,U)$ is trained on part of the
family and evaluated on simulations that were never used.

\begin{table}[!ht]
\caption{Simulations of the parametric pulsatile family. $\mu$ in Pa$\cdot$s, $U$ in m/s. Test simulations are never used for
training, scaling or checkpoint selection.}\label{tbl:param_cases}
\begin{tabular*}{\tblwidth}{@{}LLLLLL@{}}
\toprule
$\mu$ & $U=0.40$ & $U=0.55$ & $U=0.70$ & $U=0.85$ & $U=1.00$ \\
\midrule
0.003 & train & train & train & train & train \\
0.004 & train & train & \textbf{test} & train & train \\
0.005 & train & \textbf{test} & train & train & train \\
\bottomrule
\end{tabular*}
\end{table}

Thirteen of the 15 simulations are used for training. In each of them, FEBio velocities $(u,v)$ are used at a random
15\,\% of the nodes of the 2-D mesh (the thin $z$ direction of the slab is ignored) at all of that simulation's output times.
The remaining terms of the loss use no FEBio field data: the unsteady Navier--Stokes residuals at collocation points drawn
over the whole $(x,y,t,\mu,U)$ box, no-slip walls, the analytical parabolic inlet, the rest initial condition, and the RCR outlet.

\begin{figure*}[t]
    \centering
    \captionsetup[subfigure]{font=small, justification=justified, singlelinecheck=false}

    \begin{subfigure}{\textwidth}
        \centering
        \includegraphics[width=\linewidth]{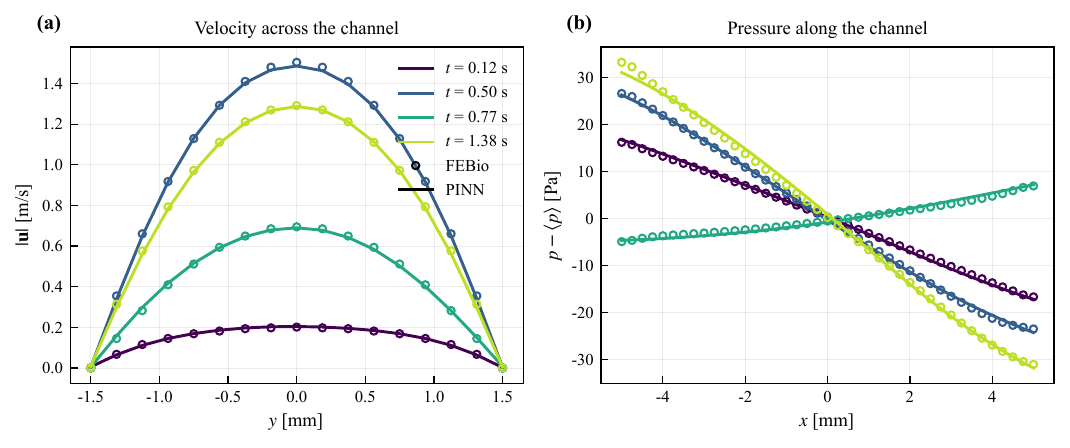}
        \caption{Velocity magnitude across the channel and pressure fluctuation along the channel at $t=0.12$, $0.50$, $0.77$
        and $1.38$~s; markers: FEBio, lines: PINN.}
        \label{fig:pulsatile_profiles}
    \end{subfigure}

    \vspace{0.3cm}

    \begin{subfigure}{\textwidth}
        \centering
        \includegraphics[width=\linewidth]{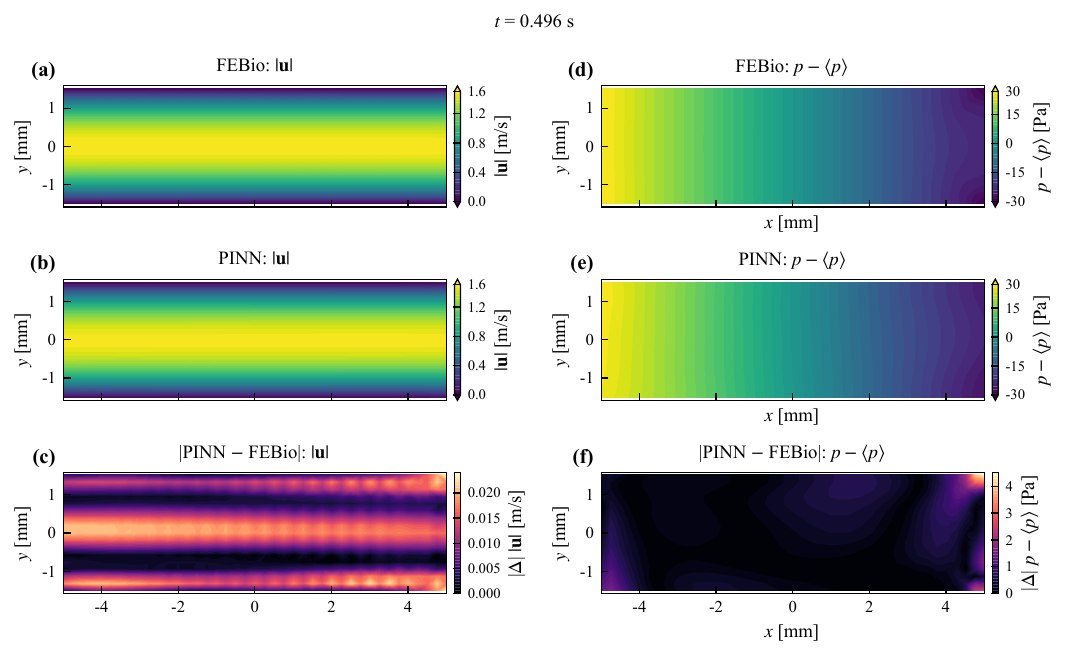}
        \caption{Fields at the instant of peak mean velocity ($t=0.496$~s). Left: velocity magnitude; right: pressure
        fluctuation $p-\langle p\rangle$. FEBio, PINN and absolute PINN--FEBio error (bottom row); FEBio and PINN share one
        colour scale per variable.}
        \label{fig:pulsatile_u}
    \end{subfigure}

    \caption{Pulsatile flow in the straight RCR channel for the unseen simulation ($\mu=0.004$~Pa$\cdot$s, $U=0.7$~m/s).}
    \label{fig:pulsatile_profiles_fields}
\end{figure*}

\begin{figure*}[t]
    \centering
    \captionsetup[subfigure]{font=small, justification=justified, singlelinecheck=false}
    \begin{subfigure}{\textwidth}
        \centering
        \includegraphics[width=\textwidth]{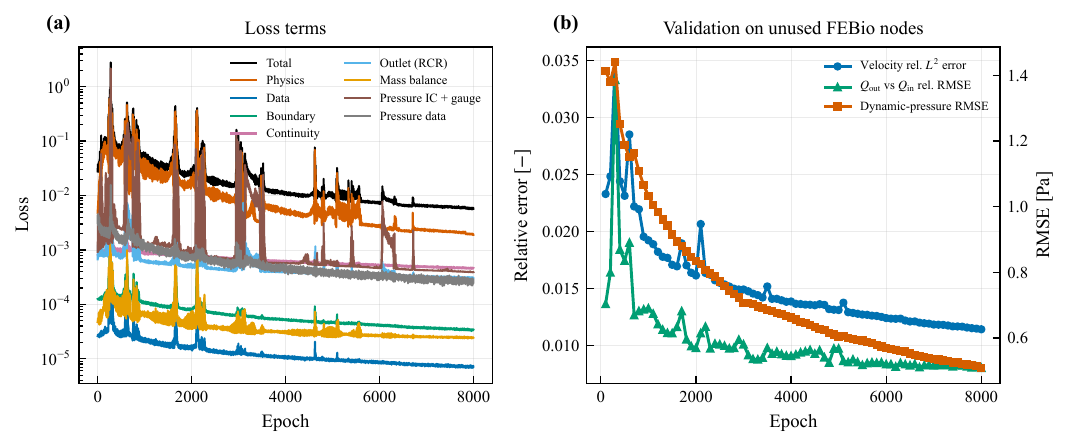}
        \caption{Terms of the loss function and validation on 400 FEBio nodes not used as data: relative $L_2$ velocity error,
        relative RMSE of $Q_{out}$ against $Q_{in}$ and RMSE of the dynamic pressure (right axis).}
        \label{fig:pulsatile_training}
    \end{subfigure}

    \vspace{0.3cm}

    \begin{subfigure}{\textwidth}
        \centering
        \includegraphics[width=\textwidth]{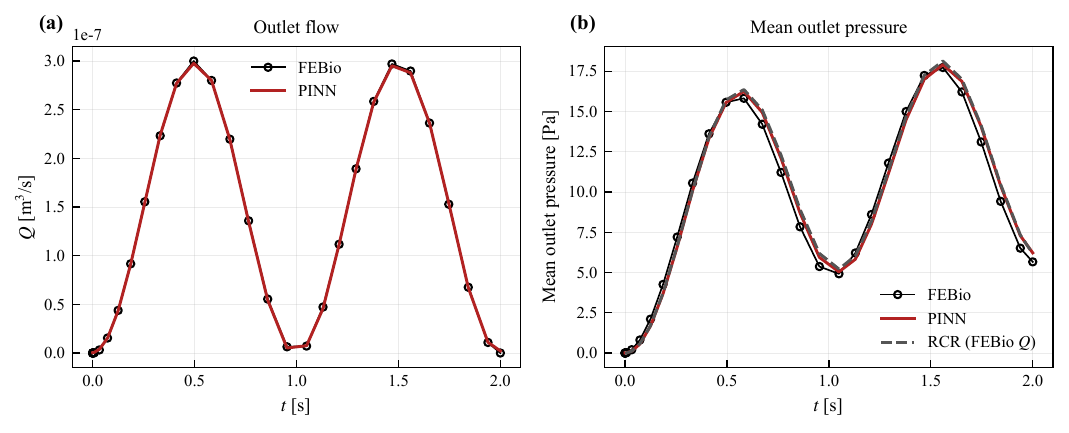}
        \caption{Outflow $Q_{out}(t)$ and mean outlet pressure; the dashed curve integrates the RCR equation with the FEBio
        flow and checks the RCR parameters and units.}
        \label{fig:pulsatile_outlet}
    \end{subfigure}

    \caption{Training of the PINN and outlet hemodynamics in the pulsatile straight channel for the unseen simulation ($\mu=0.004$~Pa$\cdot$s, $U=0.7$~m/s).}
    \label{fig:pulsatile_results}
\end{figure*}

The results below correspond to an unseen test simulation ($\mu=0.004$~Pa$\cdot$s, $U=0.7$~m/s), which was excluded from
training, scaling and checkpoint selection (Figs.~\ref{fig:pulsatile_profiles_fields} and~\ref{fig:pulsatile_results}).

Figure~\ref{fig:pulsatile_profiles} compares transverse velocity profiles (a) and axial pressure distributions (b) at four
cardiac phases ($t=0.12$, $0.50$, $0.77$ and $1.38$~s). The PINN follows the FEBio velocity through acceleration, peak systole
and deceleration. The pressure gradient has the correct sign throughout the cycle and is reproduced closely at $t=0.12$, $0.50$
and $0.77$~s; a small underestimation remains at $t=1.38$~s, where the pressure at the inlet end is about 2~Pa below the FEBio
value of $\approx33$~Pa.

Figure~\ref{fig:pulsatile_u} shows a snapshot at peak flow ($t=0.496$~s). The absolute velocity error stays below
$\approx0.02$~m/s over the whole channel, with its maximum along the bands of largest velocity gradient. The pressure error
is below $\approx1$~Pa over most of the domain and reaches $\approx4$~Pa only in the corners at the inlet and outlet ends. 

The validation history (Fig.~\ref{fig:pulsatile_training}(b)) decreases steadily and is still improving slowly at the end of
the $8000$ epochs: the relative $L_2$ velocity error is $\approx0.011$, the relative RMSE of $Q_{out}$ against $Q_{in}$ is
$\approx0.8\,\%$ and the RMSE of the dynamic pressure is $\approx0.5$~Pa. The loss terms (Fig.~\ref{fig:pulsatile_training}(a))
show transient spikes during the first $\approx3500$ epochs and are smooth afterwards, with the physics residual as the
dominant term.

Coupling with the downstream Windkessel ODE is shown in Fig.~\ref{fig:pulsatile_outlet}. The predicted outflow $Q_{out}(t)$ is
virtually identical to the FEBio solution, including the peak discharge of $3.0\times10^{-7}$~m$^3$/s in both cycles
(Fig.~\ref{fig:pulsatile_outlet}(a)); since the FE inflow enters the mass-balance term, this agreement is expected to a large
extent. The mean outlet pressure $\bar{P}_{out}(t)$ (Fig.~\ref{fig:pulsatile_outlet}(b)) follows the RCR solution driven by the
FE flow rate (dashed line) and reproduces the FEBio peaks ($\approx15.8$ and $17.7$~Pa) and trough ($\approx4.9$~Pa) to within
about $0.3$~Pa. During deceleration the nodal FE pressure deviates from the RCR solution by up to $\approx1$~Pa. We attribute
this small discrepancy to the element-to-node averaging of the FE pressure and to the time discretisation of the solver, not
to a deficiency of the outlet model.

\section{Patient-Specific Right Pulmonary Artery}\label{sec:RPA}
One of the most profound advantages of the Physics-Informed Neural Network paradigm over classical mesh-based solvers is its inherent capability to perform data assimilation and solve inverse problems seamlessly. In clinical settings, acquiring high-resolution, full-field hemodynamics is often unfeasible; instead, only sparse, noisy, or time-averaged measurements (e.g., from 4D Flow MRI or Doppler ultrasound) are typically available. 

To demonstrate the clinical translational potential of our framework, we formulate an inverse problem within the 3D patient-specific right pulmonary artery domain. In this scenario, the objective is to identify hidden physical parameters (such as the distal vascular resistance and compliance of the RCR Windkessel model), using only a sparse set of internal velocity and pressure measurements.

This section presents the application of the framework to the patient-specific three-dimensional right pulmonary artery described in Section \ref{sec:methods}. The forward problem is analyzed first in Section \ref{subsec:3d_forward}, followed by the inverse problem development in Section \ref{subsec:3d_inverse}.

\subsection{Forward problem}\label{subsec:3d_forward}
 
The reference solution was obtained with FEBio using the rigid-wall, 3WK-outlet configuration described in \cite{casero2026patient}, with a patient-specific inflow of period $1$~s simulated over two cycles. The computational domain was discretized with an unstructured mesh of 48,035 tetrahedral elements and 9,128 nodes. The boundary is split into one inflow surface $\Gamma_{in}$, the rigid wall $\Gamma_{wall}$ and one outflow surface $\Gamma_{out}$ (\texttt{FluidRCR1}) coupled to the Windkessel model. The PINN was trained with the architecture and hyper-parameters of Table~\ref{tbl:hyper}, on a workstation with an Intel Core i7-6700 CPU and an NVIDIA RTX A2000 GPU (12~GB). The $8000$ Adam epochs followed by the L-BFGS refinement took approximately 0.66~hours of wall-clock time.

\begin{figure*}[!t]
    \centering
    \includegraphics[width=\textwidth]{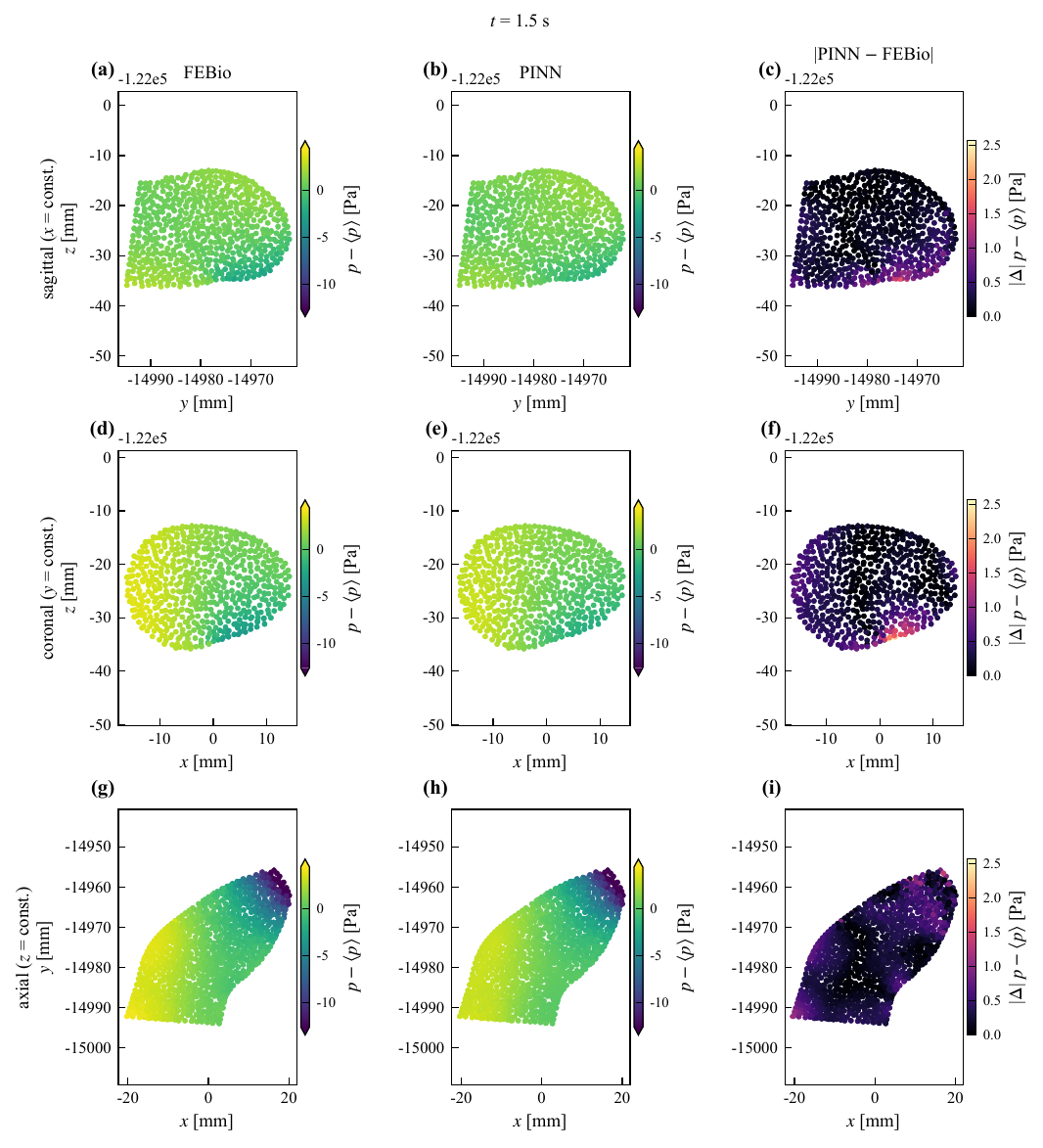}
    \caption{Dynamic pressure $p-\langle p\rangle$ in the right pulmonary artery at $t=1.4994$~s on sagittal, coronal and axial planes: FEBio (left), PINN (centre) and absolute error (right).}
    \label{fig:3d_p}
\end{figure*}

\begin{figure*}[t]
  \centering
  \includegraphics[width=\textwidth]{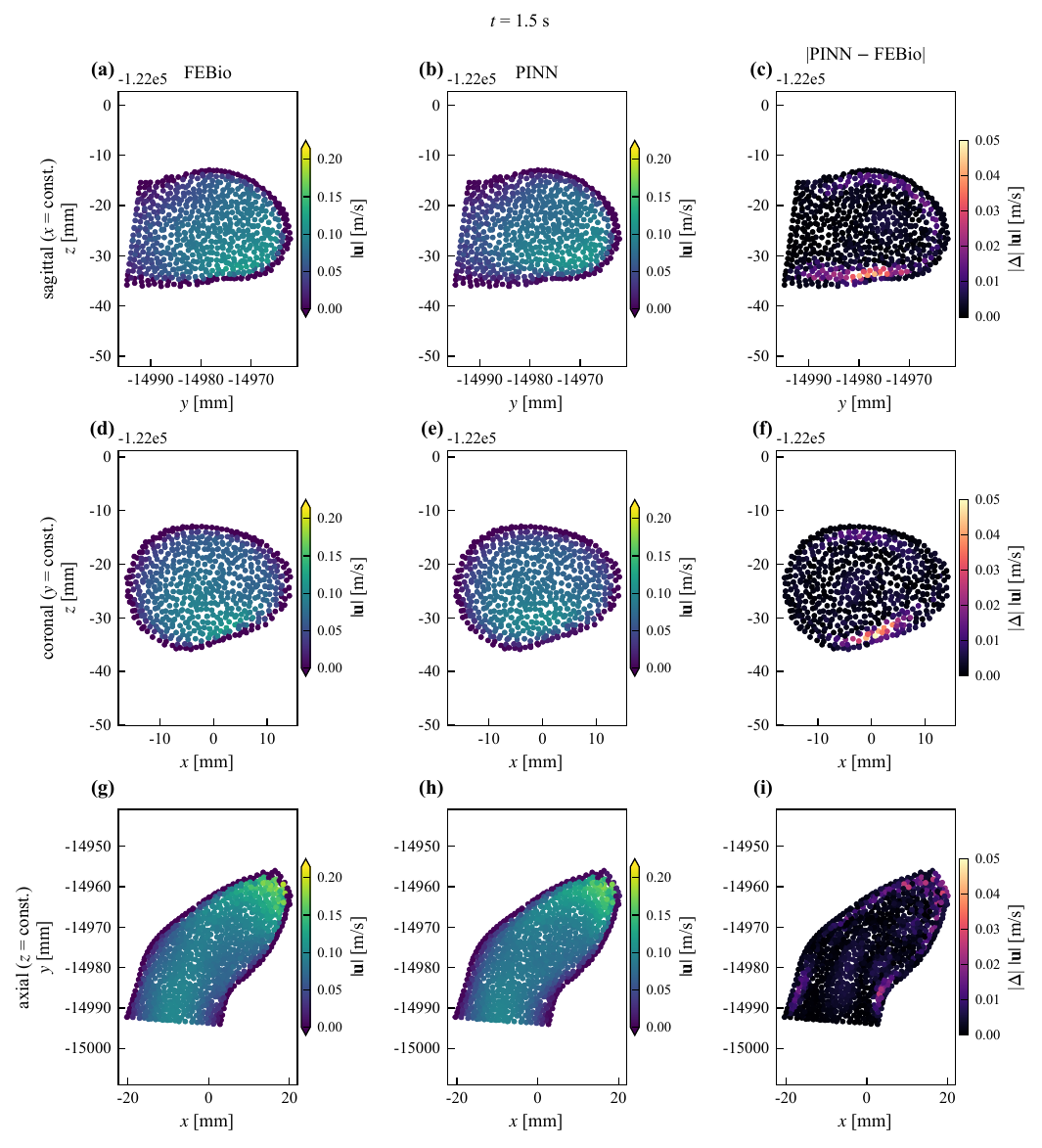}
  \caption{Velocity magnitude in the right pulmonary artery at $t=1.4994$~s on sagittal, coronal and axial planes: FEBio (left), PINN (centre) and absolute error (right).}\label{fig:3d_u}
\end{figure*}

The spatial fields at $t=1.4994$~s, close to the flow peak of the second cycle, are shown in Figs.~\ref{fig:3d_p} and \ref{fig:3d_u} for three orthogonal cutting planes (sagittal, coronal and axial). The dynamic pressure $p-\langle p\rangle$ spans approximately $[-2.3,\,2]$~Pa, $[-2.8,\,3.5]$~Pa and $[-11,\,3.5]$~Pa in the three planes, and the PINN reproduces the spatial distribution, including the pressure gradient from the proximal to the distal region. The absolute errors are below $0.3$~Pa in most of the domain and reach maxima of $1.6$, $2.5$ and $2.1$~Pa in the three planes, concentrated in the region of lowest pressure near the lower boundary. These maxima are negligible with respect to the basal pressure ($\approx0.1\%$ of $3.5$~kPa) but are not negligible with respect to the local pressure variations that drive the flow, which indicates a residual bias of the pressure field at the local scale.
 
The velocity magnitude (peak $\approx0.21$~m/s) is also correctly located and shaped (Fig.~\ref{fig:3d_u}). The velocity error is below $0.005$~m/s in the core of the vessel and reaches maxima of $0.04$, $0.05$ and $0.03$~m/s in the sagittal, coronal and axial planes, respectively, located in the near-wall layer and co-located with the regions of largest pressure error. This is the expected behavior for a network that must resolve the steep near-wall velocity gradients imposed by the no-slip condition from scattered points.
 
\begin{figure*}[t]
    \centering
    \captionsetup[subfigure]{font=small, justification=justified, singlelinecheck=false}

    \begin{subfigure}{\textwidth}
        \centering
        \includegraphics[width=\textwidth]{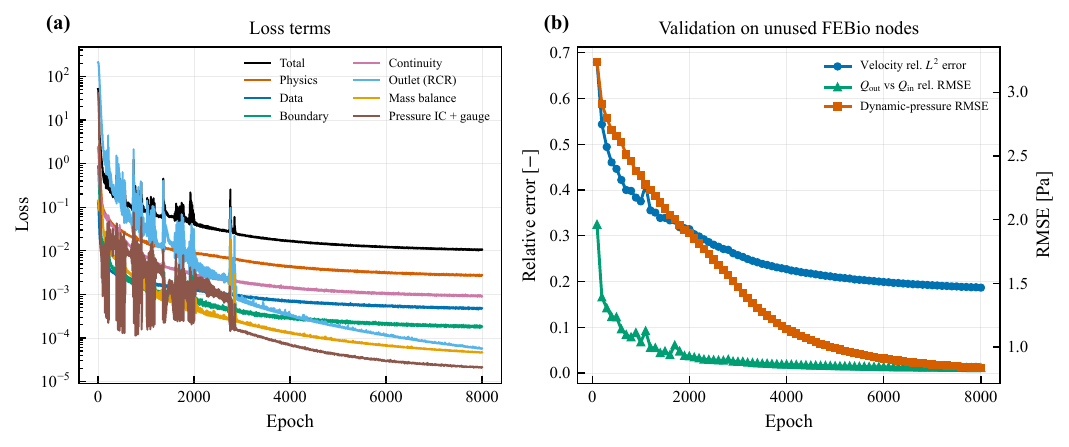}
        \caption{Training history for the right pulmonary artery model. Left: loss terms. Right: validation on held-out FE
        nodes (relative $L_2$ velocity error, relative RMSE between $Q_{out}$ and $Q_{in}$, and dynamic-pressure RMSE).}
        \label{fig:3d_training}
    \end{subfigure}

    \vspace{0.3cm}

    \begin{subfigure}{\textwidth}
        \centering
        \includegraphics[width=\textwidth]{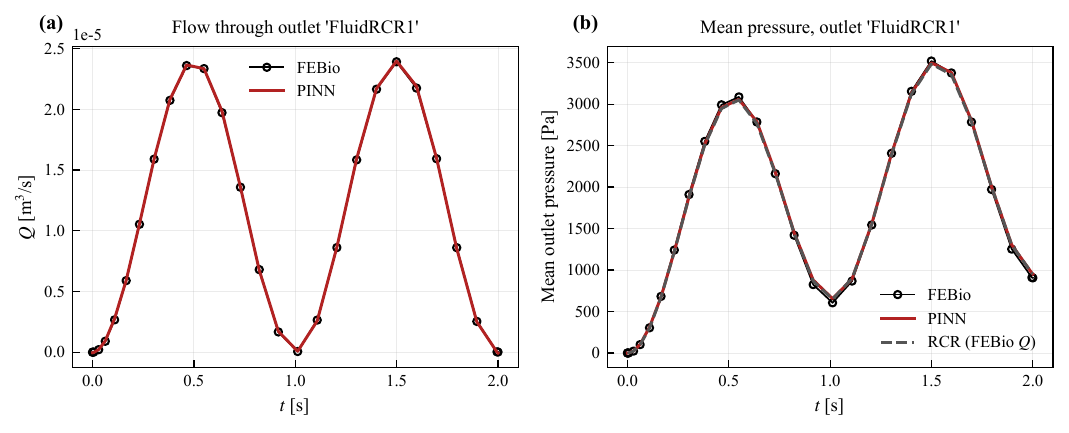}
        \caption{Outlet hemodynamics in the right pulmonary artery model. Left: flow rate through outlet
        \texttt{FluidRCR1}. Right: area-averaged outlet pressure from FEBio, from the PINN, and from the RCR ordinary
        differential equation driven by the FE flow.}
        \label{fig:3d_outlet}
    \end{subfigure}

    \caption{Forward problem in the three-dimensional right pulmonary artery model: training history and outlet hemodynamics.}
    \label{fig:3d_results}
\end{figure*}

 The training history (Fig.~\ref{fig:3d_training}) shows a monotonic decrease of all loss terms after the initial transient, which ends once the stiff terms are fully active ($E_r=1500$ epochs) and, for the Windkessel and anchor terms, around $2800$ epochs. At the end of training, the relative $L_2$ velocity error on the held-out nodes is $\approx0.19$, the dynamic-pressure RMSE is $\approx0.83$~Pa and the relative RMSE between $\sum_kQ_{out}^{(k)}$ and $Q_{in}$ is $\approx1\%$. The velocity error is thus about one order of magnitude larger than in the two-dimensional pulsatile case ($0.02$), which reflects the increased difficulty of the three-dimensional geometry, of the near-wall resolution and of the limited fraction of supervised nodes ($15\%$). 

 Figure~\ref{fig:3d_outlet} also shows the flow rate through the outlet and the area-averaged outlet pressure. The outflow is reproduced essentially exactly (peak $\approx2.4\times10^{-5}$~m$^3$/s in both cycles). The mean outlet pressure reaches $\approx3.1$~kPa in the first cycle and $\approx3.5$~kPa in the second, and the PINN matches the FE value within $\approx1.5\%$ at the peaks and within $\approx50$~Pa at the diastolic trough. The pressure also coincides with the solution of the RCR equation driven by the FE flow rate, which shows that the network satisfies the Windkessel relation consistently. The higher pressure of the second cycle with respect to the first is a consequence of the charging of the compliance from rest, which is represented by the time-only basal-pressure network (Section \ref{subsec:pinn_arch}). The Windkessel-imposed absolute pressure is therefore three orders of magnitude larger than the local pressure variations that drive the flow, which are of the order of $1$--$10$~Pa; this is precisely the scale disparity that the pressure decomposition of Eq. \eqref{eq:pdecomp} is designed to handle.

The results show that the proposed decomposition of the pressure and the dedicated treatment of the Windkessel relation allow the network to recover both the global (outlet) hemodynamics and the local fields in a realistic geometry, despite the three-order-of-magnitude difference between basal and dynamic pressures. The remaining limitations are the residual bias of the local pressure (errors up to $2.5$~Pa compared with a dynamic range of a few pascals) and the near-wall velocity error. Both are consistent with the behavior observed in the canonical cases: the pressure is the quantity least constrained by the data and the physics, and it is the last to converge.

\clearpage
\subsection{Inverse problem}\label{subsec:3d_inverse}

The inverse problem is formulated using synthetic, noise-free observations extracted from the finite element (FE) solution of the forward problem. Three types of measurements are assumed to be available: (i) velocity at the same sparse subset of FE nodes considered in the forward problem, corresponding to $15\%$ of the total nodes and available at all output times; (ii) the volumetric outflow rate, $Q_{out}(t)$, computed from the velocity field at the outlet; and (iii) the area-averaged outlet pressure, $\bar P_{out}(t)$, representing measurements that could be obtained from an implantable pressure sensor. No pressure measurements are prescribed within the computational domain, and the weight associated with the pressure data term is therefore set to zero ($\lambda_p=0$).

The parameters to be identified are the proximal resistance $R_p$, distal resistance $R_d$, and compliance $C$ of the outlet \texttt{FluidRCR1} Windkessel model. The parameter-identification problem is formulated by promoting these quantities to trainable variables, which are optimized jointly with the neural network weights $\theta$. The corresponding Windkessel relation is enforced at the outlet time instants $t_0<t_1<\dots<t_N$ using a backward-Euler discretization of the time derivatives:

\begin{equation}\label{eq:rcr_be}
\begin{aligned}
    r_k ={}& C\,\frac{\bar P_k-\bar P_{k-1}}{\Delta t_k}
    +\frac{\bar P_k-P_d}{R_d} \\
    &-Q_k\Big(1+\frac{R_p}{R_d}\Big)
    -R_pC\,\frac{Q_k-Q_{k-1}}{\Delta t_k},
\end{aligned}
\end{equation}

where $\Delta t_k=t_k-t_{k-1}$, and $\bar P_k$ and $Q_k$ denote the area-averaged pressure and volumetric flow rate predicted by the network at time $t_k$, respectively (Eq.~\eqref{eq:Qdisc}). The corresponding residual loss is defined as

\begin{equation}
    \mathcal{L}_{rcr} =
    \frac{\langle r_k^2\rangle}{Q_{ref}^2}.
\end{equation}

\begin{figure*}[htbp]
    \centering
    \includegraphics[width=\linewidth]{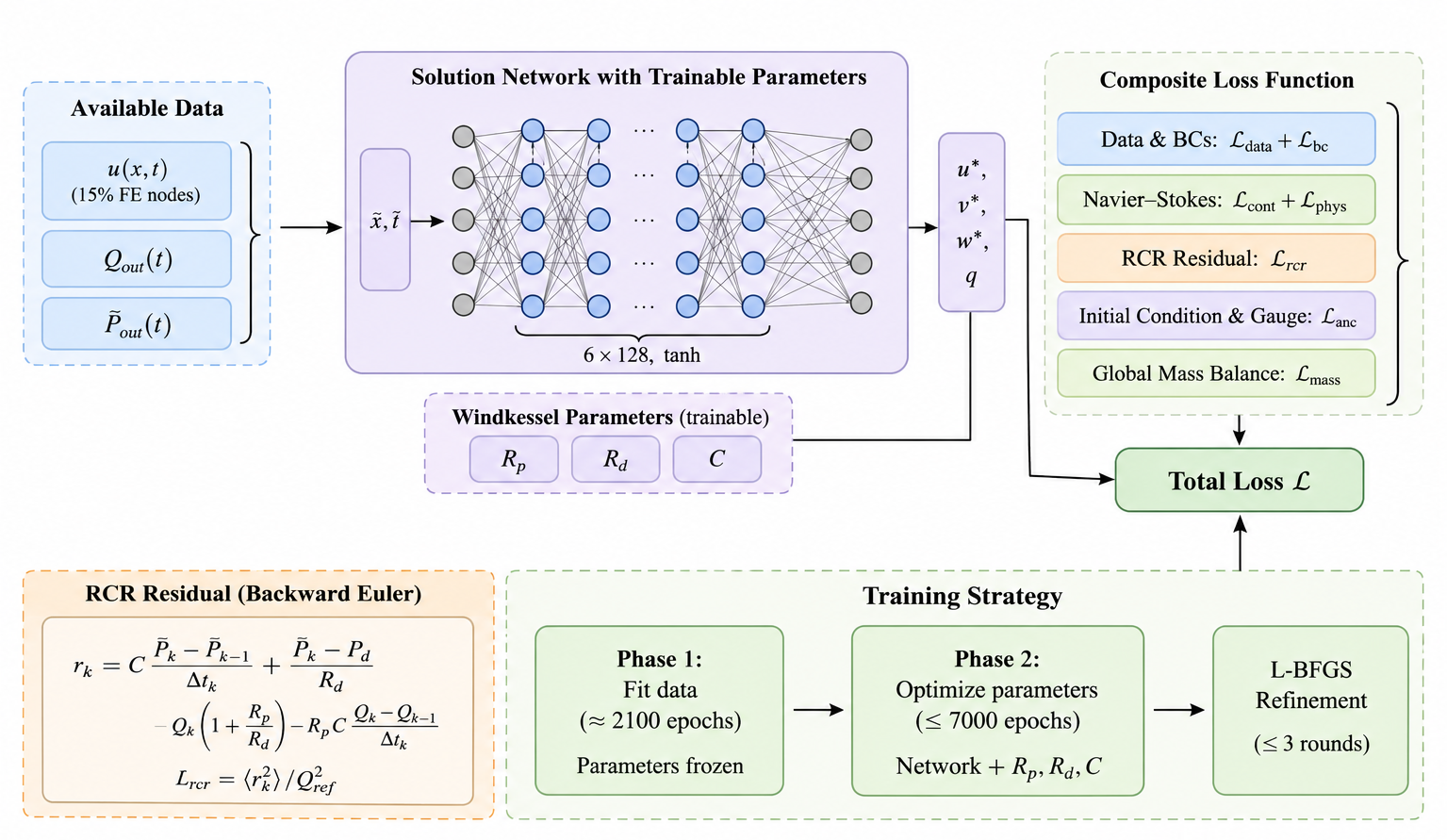}
    \caption{Schematic representation of the DAPINN framework for the inverse problem. The Windkessel parameters $(R_p,R_d,C)$ are treated as trainable variables and identified jointly with the neural network weights. The RCR relation is enforced through a backward-Euler discretization of the outlet model.}
    \label{fig:DAPINN}
\end{figure*}

The optimization is performed in two stages. During the first approximately $2100$ epochs, the Windkessel parameters are held fixed at deliberately perturbed initial values, allowing the neural network to first adapt to the available observations. Subsequently, the parameters are released and jointly optimized with the network weights until epoch $7000$. This is followed by additional refinement using L-BFGS optimization. The overall DAPINN framework is illustrated in Fig.~\ref{fig:DAPINN}.

The reference parameter values, obtained from the FE model, are $R_p=1.1\times10^{8}$~Pa$\cdot$s/m$^3$, $R_d=5.21\times10^{7}$~Pa$\cdot$s/m$^3$, and $C=1.28\times10^{-8}$~m$^3$/Pa. To assess the robustness of the identification procedure, additional test cases are generated by modifying one or more of the reference Windkessel parameters. These cases are labelled \texttt{AP1\_<parameter>\_x<factor>}, where the indicated parameter(s) are multiplied by the specified factor in the FE model, thereby generating a new synthetic dataset with a corresponding ground truth. For example, \texttt{AP1\_Rp\_x0.5} corresponds to a simulation with $R_p=0.5R_p^{\mathrm{ref}}$.

As an independent reference, the identified parameters are also compared against those obtained through a reduced nonlinear least-squares (NLS) identification of the RCR model using the same outlet flow and pressure measurements.

\begin{figure*}[t]
  \centering
  \includegraphics[width=\textwidth]{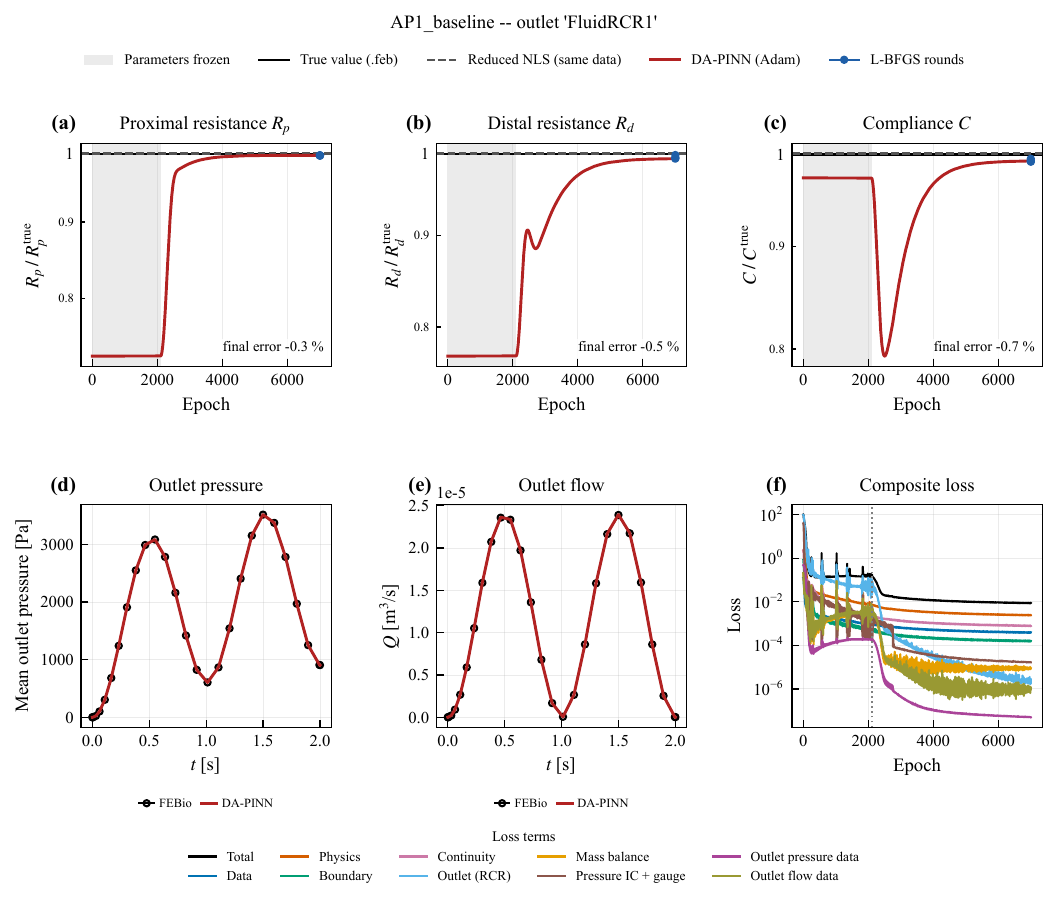}
  \caption{Inverse problem for the outlet Windkessel parameters. Top: evolution of the estimated $R_p$, $R_d$, and $C$, normalized by their respective ground-truth values (grey band: epochs during which the parameters are held fixed; red markers: L-BFGS refinement rounds; dashed lines: parameter estimates obtained through the reduced NLS fit using the same observations). Bottom: reconstructed mean outlet pressure, outlet flow rate, and loss contributions.}
  \label{fig:inverse}
\end{figure*}

The evolution of the identified parameters for the reference case is presented in Fig.~\ref{fig:inverse}. At initialization, $R_p$ is underestimated by approximately $34\%$, whereas $R_d$ is overestimated by approximately $40\%$. Following the release of the Windkessel parameters, both resistances converge towards their reference values within approximately $1000$ epochs. The compliance exhibits a slower convergence process, initially decreasing below its reference value before progressively approaching it towards the end of training.

The final estimates, reported in Table~\ref{tbl:inverse}, show that all three parameters are recovered with relative errors below $0.3\%$. The reconstructed outlet pressure and flow waveforms closely reproduce the corresponding synthetic observations (Fig.~\ref{fig:inverse}, bottom left and centre). Moreover, releasing the Windkessel parameters leads to an approximately one-order-of-magnitude reduction in the total loss, while the outlet-related loss contributions decrease to the order of $10^{-6}$ (Fig.~\ref{fig:inverse}, bottom right).

The different convergence rates observed for the three parameters suggest that the compliance is less readily identified than the resistances from the available pressure and flow waveforms. This behaviour is consistent with the different sensitivities of the RCR model outputs to its individual parameters and may become particularly relevant when considering measurement noise and parameter identifiability.

\begin{table}[htbp]
\caption{Identified Windkessel parameters for the reference inverse problem using noise-free synthetic data.}
\label{tbl:inverse}
\begin{tabular*}{\tblwidth}{@{}LLLL@{}}
\toprule
Parameter & Ground truth & Identified & Relative error \\
\midrule
$R_p$ [Pa$\cdot$s/m$^3$] & $1.10\times10^{8}$ & $1.097\times10^{8}$ & $-0.3\%$ \\
$R_d$ [Pa$\cdot$s/m$^3$] & $5.21\times10^{7}$ & $5.22\times10^{7}$ & $+0.2\%$ \\
$C$ [m$^3$/Pa] & $1.28\times10^{-8}$ & $1.276\times10^{-8}$ & $-0.3\%$ \\
\bottomrule
\end{tabular*}
\end{table}

\begin{figure*}[!t]
  \centering
  \subfloat[For $R_p\times0.5$\label{fig:inverse_Rp05}]{%
  \includegraphics[width=0.9\textwidth, trim={0 1.2cm 0 1.8cm}, clip]{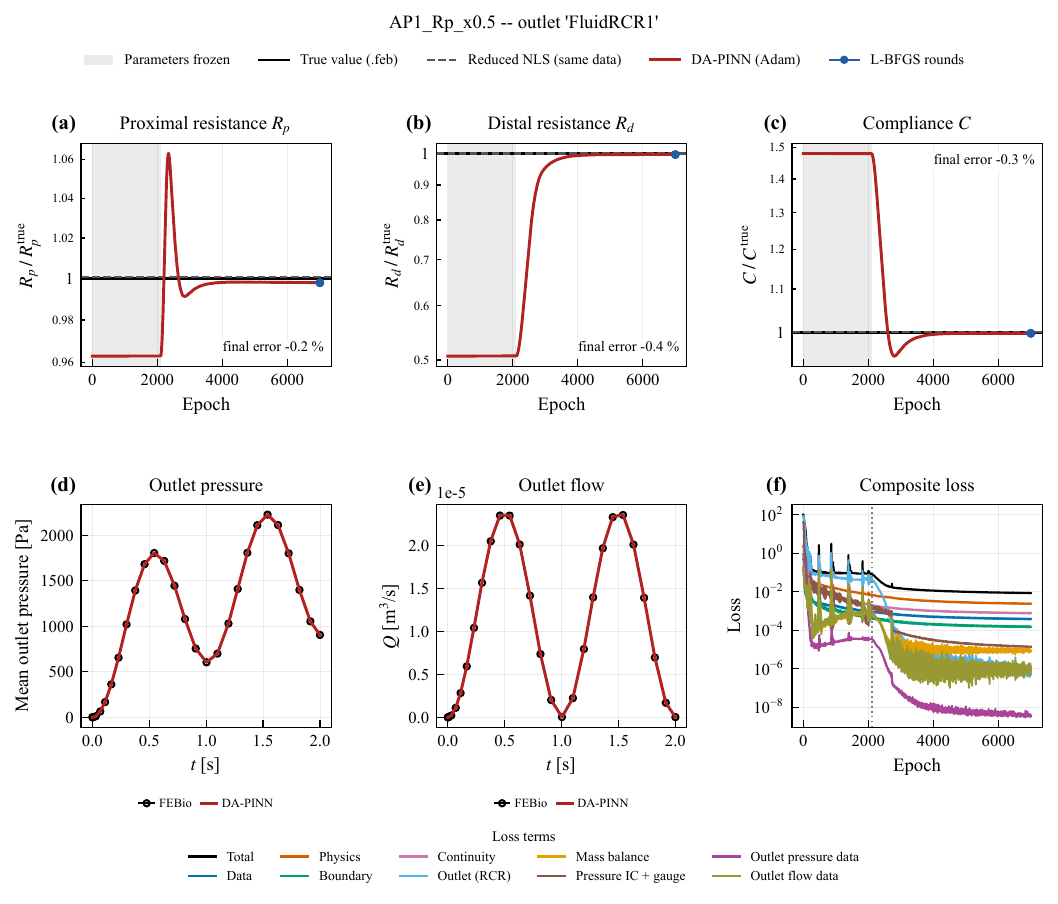}}\\[0.15cm]
  \subfloat[For $R_d\times2.0$\label{fig:inverse_Rd2}]{%
  \includegraphics[width=0.9\textwidth, trim={0 1.2cm 0 1.8cm}, clip]{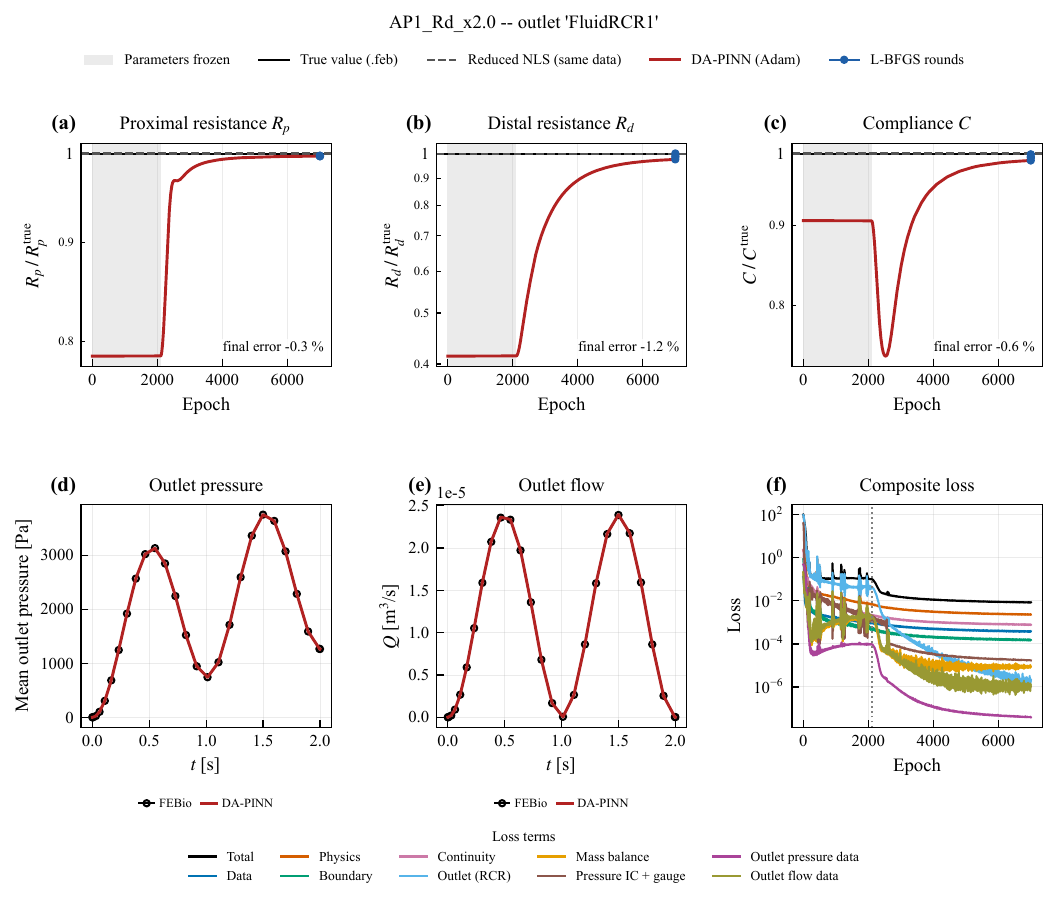}}
  \caption{Parameter identification results for modified ground-truth Windkessel parameters: (a) $R_p\times0.5$ and (b) $R_d\times2.0$. The layout follows that of Fig.~\ref{fig:inverse}.}
  \label{fig:inverse_variants_a}
\end{figure*}

\begin{figure*}[!t]
  \centering
  \subfloat[For $C\times0.5$\label{fig:inverse_C05}]{%
    \includegraphics[width=0.9\textwidth, trim={0 1.2cm 0 1.8cm}, clip]{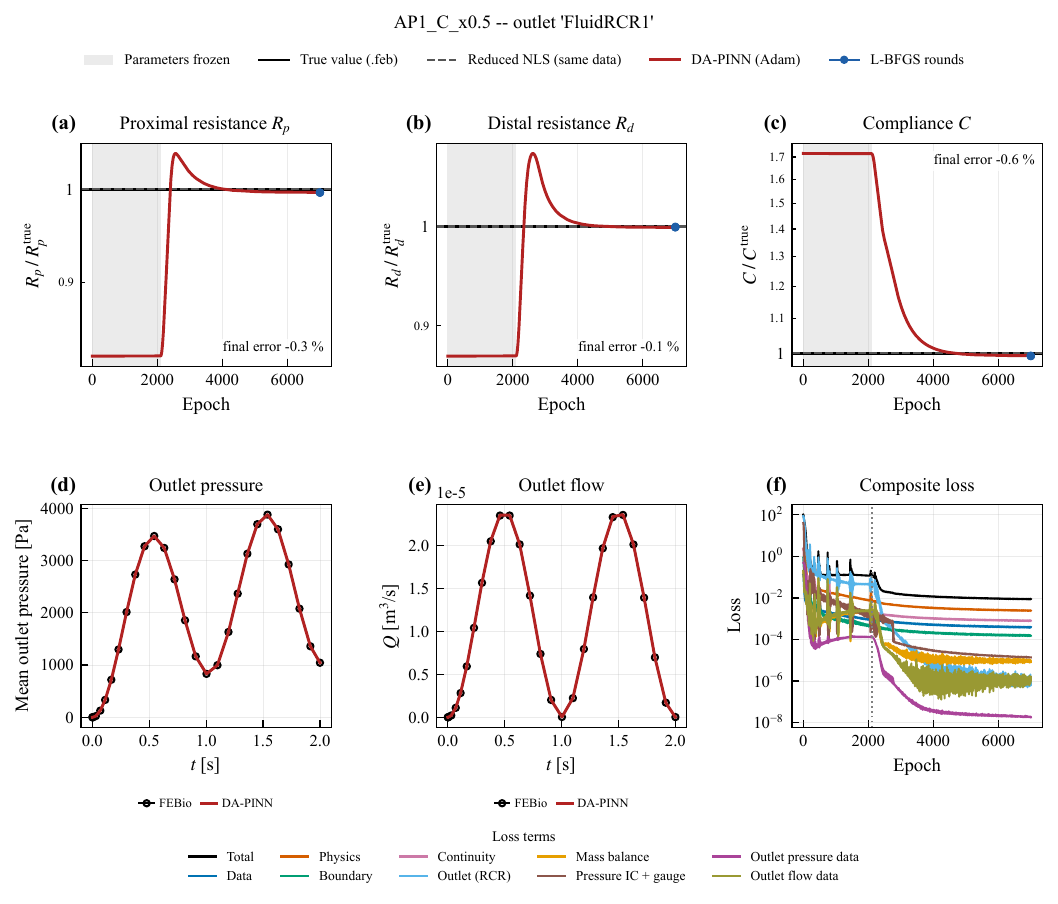}}\\[0.15cm]
  \subfloat[For $R_p\times2.0$ and $C\times0.5$\label{fig:inverse_Rp2C05}]{%
    \includegraphics[width=0.9\textwidth, trim={0 1.2cm 0 1.8cm}, clip]{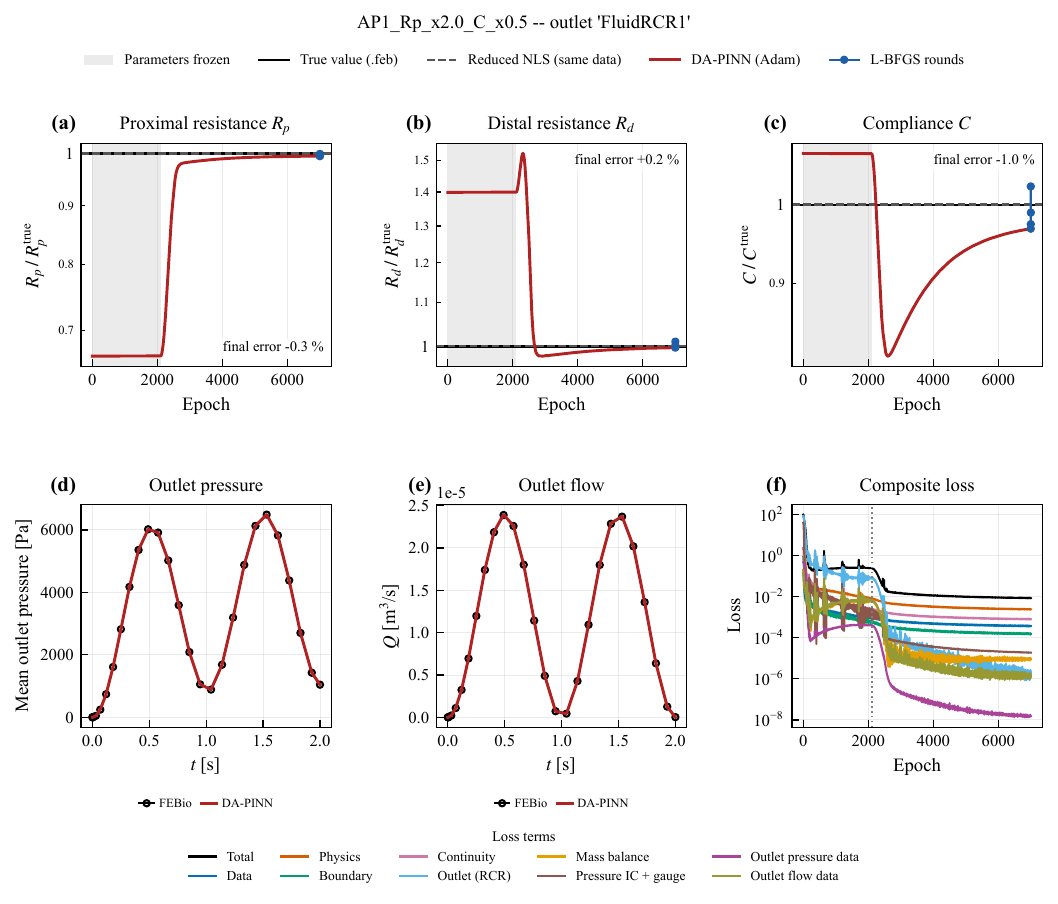}}
  \caption{Parameter identification results for modified ground-truth Windkessel parameters: (a) $C\times0.5$ and (b) the combined variation $R_p\times2.0$, $C\times0.5$. The layout follows that of Fig.~\ref{fig:inverse}.}
  \label{fig:inverse_variants_b}
\end{figure*}

The performance of the identification framework under different parameter configurations is further investigated in Figs.~\ref{fig:inverse_variants_a} and~\ref{fig:inverse_variants_b}, with the corresponding quantitative results summarized in Table~\ref{tbl:inverse_variants}. Four additional scenarios are considered: $R_p\times0.5$, $R_d\times2.0$, $C\times0.5$, and the combined variation $R_p\times2.0$ and $C\times0.5$.

These modifications lead to substantial changes in the outlet pressure response, with peak mean pressures ranging from approximately $2.2$~kPa for the $R_p\times0.5$ case to $6.5$~kPa for the combined $R_p\times2.0$, $C\times0.5$ case. In contrast, the outlet flow remains practically unchanged, with a peak value of approximately $2.4\times10^{-5}$~m$^3$/s, consistent with the prescribed inflow conditions.

Across all four scenarios, the framework successfully recovers the three Windkessel parameters, with final relative errors ranging from $0.1\%$ to $1.2\%$ and a mean absolute relative error of approximately $0.5\%$. These results are obtained despite substantial deviations in the initial parameter estimates, which reach approximately $58\%$ for $R_d$ in the $R_d\times2.0$ case and $71\%$ for $C$ in the $C\times0.5$ case. In all scenarios, the reconstructed outlet pressure and flow waveforms closely match the corresponding observations. Furthermore, following the release of the parameters, the outlet-related loss contributions decrease to approximately $10^{-6}$, with the mean pressure and flow misfit terms reaching values between $10^{-8}$ and $10^{-6}$.

Several consistent trends can be identified across the different parameter configurations. First, $R_p$ exhibits a small systematic negative bias in all cases, with relative errors between $-0.2\%$ and $-0.3\%$. This consistency suggests that the bias may be associated with the numerical formulation of the discretized RCR relation, potentially arising from the first-order backward-Euler approximation, rather than being solely attributable to optimization variability.

Second, $R_d$ and $C$ exhibit slower convergence and more pronounced transient deviations following the release of the parameters. In particular, transient deviations of up to approximately $40\%$ are observed for $C$ in the $R_d\times2.0$ and combined $R_p\times2.0$, $C\times0.5$ scenarios. This behaviour suggests a degree of compensation between the resistance and compliance parameters during optimization, reflecting their coupled influence on the outlet pressure and flow responses.

Finally, in the scenarios involving the largest initial parameter deviations, namely $R_d\times2.0$ and the combined $R_p\times2.0$, $C\times0.5$ case, the estimates of $R_d$ and $C$ are still evolving towards their respective ground-truth values at epoch $7000$. Consequently, the remaining errors of up to $1.2\%$ should be interpreted in the context of incomplete convergence at the prescribed training horizon. Further optimization may therefore lead to additional improvements in the identified values.

Overall, these results demonstrate the ability of the proposed DAPINN framework to identify the lumped Windkessel parameters from sparse velocity observations and outlet pressure and flow measurements, even when the initial parameter estimates differ substantially from the ground truth. Importantly, the identification is achieved through joint optimization of the neural network and the outlet model parameters, without requiring repeated FE forward simulations during the inverse optimization.

\begin{table*}[htbp]
\caption{Identification of the Windkessel parameters for the four modified scenarios using noise-free synthetic data (7000 epochs followed by L-BFGS refinement). Relative errors are calculated as $(\hat{\theta}-\theta)/\theta$. Initial errors represent approximate deviations of the data-driven initial estimates from the ground-truth values, as inferred from the parameter evolution plots.}
\label{tbl:inverse_variants}
\small
\begin{tabular*}{\textwidth}{@{\extracolsep{\fill}}lccc ccc ccc@{}}
\toprule
 & \multicolumn{3}{c}{Ground truth} & \multicolumn{3}{c}{Initial error (approx.)} & \multicolumn{3}{c}{Final relative error} \\
\cmidrule(lr){2-4}\cmidrule(lr){5-7}\cmidrule(lr){8-10}
Case & $R_p$ [Pa$\cdot$s/m$^3$] & $R_d$ [Pa$\cdot$s/m$^3$] & $C$ [m$^3$/Pa] & $R_p$ & $R_d$ & $C$ & $R_p$ & $R_d$ & $C$ \\
\midrule
$R_p\times0.5$ & $5.48\times10^{7}$ & $1.04\times10^{8}$ & $1.28\times10^{-8}$ & $-4\%$ & $-49\%$ & $+48\%$ & $-0.2\%$ & $-0.4\%$ & $-0.3\%$ \\
$R_d\times2.0$ & $1.10\times10^{8}$ & $2.08\times10^{8}$ & $1.28\times10^{-8}$ & $-21\%$ & $-58\%$ & $-10\%$ & $-0.3\%$ & $-1.2\%$ & $-0.6\%$ \\
$C\times0.5$ & $1.10\times10^{8}$ & $1.04\times10^{8}$ & $6.42\times10^{-9}$ & $-17\%$ & $-13\%$ & $+71\%$ & $-0.3\%$ & $-0.1\%$ & $-0.6\%$ \\
$R_p\times2.0,\ C\times0.5$ & $2.19\times10^{8}$ & $1.04\times10^{8}$ & $6.40\times10^{-9}$ & $-34\%$ & $+40\%$ & $+7\%$ & $-0.3\%$ & $+0.2\%$ & $-1.0\%$ \\
\bottomrule
\end{tabular*}
\end{table*}

\FloatBarrier
\clearpage
\section{Conclusions and Future Perspectives}\label{sec:conclusions}

This work presented a physics-informed neural network (PINN) framework for modelling pulsatile blood flow in patient-specific pulmonary artery geometries, with particular emphasis on the integration of dynamic three-element Windkessel (RCR) outlet boundary conditions. The proposed approach addresses the numerical challenges associated with coupling local three-dimensional hemodynamics to lumped-parameter models, particularly the disparity between the pressure scales governing the intravascular flow and those associated with the outlet boundary conditions. A progressive validation strategy was adopted, ranging from canonical two-dimensional configurations to a patient-specific three-dimensional right pulmonary artery, followed by an inverse problem for the identification of the outlet model parameters.

The initial validation on canonical geometries demonstrated the ability of the framework to reproduce both analytical solutions and FEBio reference simulations under steady and transient flow conditions. For the straight-channel benchmark, the predicted velocity profiles agreed closely with the analytical Poiseuille solution and the corresponding FE results, with maximum pointwise velocity errors below $0.6\%$ of the peak velocity. The extension to the channel with lateral offset demonstrated that the network could also capture more complex flow distributions, including local acceleration, curved streamlines and non-uniform pressure gradients. These tests additionally highlighted the limitations of relying exclusively on a soft incompressibility constraint to preserve global mass conservation, motivating the incorporation of an explicit mass-balance loss term. Finally, the pulsatile benchmark confirmed that the framework could reproduce transient velocity and pressure distributions while consistently coupling the network predictions to the downstream RCR model.

The application to the patient-specific three-dimensional right pulmonary artery further demonstrated the feasibility of the proposed approach in a geometrically complex domain. The network reproduced the main spatial features of the FE velocity and dynamic pressure fields, as well as the temporal evolution of the outlet flow and mean pressure. In particular, the reconstructed outlet pressure was consistent with the solution of the RCR ordinary differential equation driven by the FE flow rate. The complete three-dimensional simulation was trained in approximately $0.66$ hours, illustrating the computational feasibility of the proposed formulation for this application.

Beyond the forward problem, the framework was extended to the identification of the proximal resistance, distal resistance and compliance of the outlet Windkessel model. Using sparse velocity observations, volumetric outflow and area-averaged outlet pressure, the three parameters were optimized jointly with the neural network weights, without requiring repeated FE simulations during the inverse optimization. For the reference configuration, all parameters were recovered with relative errors below $0.3\%$. Additional experiments involving substantial modifications of the ground-truth parameters demonstrated consistent parameter recovery, with final errors ranging from $0.1\%$ to $1.2\%$ across the four investigated scenarios. These results demonstrate the potential of the proposed formulation to combine distributed hemodynamic reconstruction with the identification of hidden lumped-parameter boundary conditions.

Despite these encouraging results, several limitations remain to be addressed. First, the inverse problem was evaluated using synthetic, noise-free observations generated from the FE model. The sensitivity of parameter identification to measurement noise, temporal undersampling and uncertainty in the available observations therefore remains to be investigated. Second, although the overall pressure distributions were satisfactorily reconstructed, localized discrepancies persisted in regions with strong spatial gradients and near geometric discontinuities. Further improvements in pressure reconstruction and parameter identifiability will be particularly important when only sparse clinical measurements are available. Finally, the present study assumes rigid vessel walls and a prescribed inflow condition, limiting its representation of the dynamic interaction between blood flow and the surrounding arterial tissue.

Future research will focus on advancing the proposed pipeline towards clinical translation. Key directions include the assimilation of sparse clinical velocity measurements, such as those obtained from 4D Flow MRI or Doppler ultrasound, directly into the PINN loss function; the incorporation of arterial wall compliance through physics-informed fluid--structure interaction (FSI) formulations; and the use of neural operators, such as DeepONets and Fourier Neural Operators, to further accelerate model training and inference. In particular, combining these developments with the inverse identification capabilities presented in this work could enable the integration of patient-specific anatomical, hemodynamic and vascular mechanical information within a unified computational framework.

Overall, the results establish a foundation for integrating physics-informed learning with dynamic lumped-parameter cardiovascular models in three-dimensional patient-specific domains. By combining forward hemodynamic reconstruction and inverse parameter identification within a single optimization framework, the proposed approach represents a promising step towards computationally efficient, non-invasive and patient-specific cardiovascular digital twins. Such developments could ultimately support real-time hemodynamic assessment and personalized modelling in pulmonary vascular medicine.


\clearpage
\printcredits

\section*{Ethics approval}
This study was conducted in accordance with the Declaration of Helsinki, and the protocol was approved by the Research Ethics Committee of the University Hospital Puerta de Hierro Majadahonda on 1 February 2024. Part of the data was collected within the framework of the academic internships of Laura García Sánchez (September--December 2022) and Angélica Casero Baeza (September--December 2023) at the Hospital Universitario Puerta de Hierro Majadahonda. Informed consent was obtained from the subject exclusively for the use of anonymized medical imaging and inflow data. In this research, synthetic pressure data was employed.

\section*{Declaration of Generative AI and AI-assisted technologies in the writing process}
During the preparation of this work, the authors used Gemini in order to reorganize paragraphs and correct spelling and grammar. After using this tool/service, the authors reviewed and edited the content as needed and take full responsibility for the content of the published article.

\section*{Funding}
P. Navas, C. Taboada, and A. Casero were supported by DOCTORES-EMERGENTES-24-6R4VJE-38-PAQCQO (TuCoPINN).

\section*{Declaration of competing interest}
The authors declare that they have no known competing financial interests or personal relationships that could have appeared to influence the work reported in this paper.

\section*{Data availability}
Data and codes will be made available upon reasonable request.

\section*{Acknowledgements}
The authors would like to express their gratitude to engineer Laura García Sánchez for her foundational role in initiating this work. We also extend our sincere thanks to Dr. Juan Francisco Oteo, Dr. Carlos Arellano Serrano, and Dr. Manuel Gómez Bueno for their invaluable clinical insights and for their active participation in the project by hosting and mentoring students at the Hospital Universitario Puerta de Hierro Majadahonda.

\bibliographystyle{cas-model2-names}

\bibliography{hemo-refs}

\end{document}